\documentclass[12pt,a4paper]{amsart}

\usepackage{amsfonts,amsmath,amssymb}
\def\frak{\mathfrak}

\newtheorem{theorem}{Theorem}[section]
\newtheorem{lemma}[theorem]{Lemma}
\newtheorem{corollary}[theorem]{Corollary}
\newtheorem{conjecture}[theorem]{Conjecture}
\newtheorem{proposition}[theorem]{Proposition}

\theoremstyle{definition}

\newtheorem{example}[theorem]{Example}

\theoremstyle{remark}
\newtheorem{remark}[theorem]{Remark}

\numberwithin{equation}{section}

\def\br{\mathbb R}
\def\bc{\mathbb C}
\def\bh{\mathbb H}
\def\bz{\mathbb Z}

\def\mand{\quad\mbox{and}\quad}
\def\mbx#1{\makebox[#1cm]{}}

\def\Ad{{\rm Ad}}
\def\ad{{\rm ad}}
\def\bp{\begin{pmatrix}}
\def\ep{\end{pmatrix}}
\def\op{{\rm op}}

\begin{document}
\bibliographystyle{alpha}

\title[Equivalence of domains]{Equivalence of domains arising from duality of orbits on flag manifolds}
\author{Toshihiko MATSUKI}
\address{\hskip-\parindent
        Toshihiko Matsuki\\
        Department of Mathematics\\
        Faculty of Science\\
        Kyoto University\\
        %Street address or POBox\\
        Kyoto 606-8502, Japan}
\email{matsuki@math.kyoto-u.ac.jp}
\date{}

\begin{abstract}
In \cite{GM1}, we defined a $G_\br$-$K_\bc$ invariant subset $C(S)$ of $G_\bc$ for each $K_\bc$-orbit $S$ on every flag manifold $G_\bc/P$ and conjectured that the connected component $C(S)_0$ of the identity would be equal to the Akhiezer-Gindikin domain $D$ if $S$ is of nonholomorphic type by computing many examples. In this paper, we first prove (in Theorem 1.3 and Corollary 1.4) this conjecture for the open $K_\bc$-orbit $S$ on an ``arbitrary'' flag manifold generalizing the result of Barchini. This conjecture for closed $S$ was solved in \cite{WZ1}, \cite{WZ2} (Hermitian cases) and \cite{FH} (non-Hermitian cases). We also deduce an alternative proof of this result for non-Hermitian cases from Theorem 1.3.
\end{abstract}

\maketitle

\section{Introduction}

Let $G_\bc$ be a connected complex semisimple Lie group and $G_\br$ a connected real form of $G_\bc$. Let $K_\bc$ be the complexification in $G_\bc$ of a maximal compact subgroup $K$ of $G_\br$. Let $X=G_\bc/P$ be a flag manifold of $G_\bc$ where $P$ is an arbitrary parabolic subgroup of $G_\bc$. Then there exists a natural one-to-one correspondence between the set of $K_\bc$-orbits $S$ and the set of $G_\br$-orbits $S'$ on $X$ given by the condition:
\begin{equation}
S\leftrightarrow S'\Longleftrightarrow S\cap S'\mbox{ is non-empty and compact} \tag{1.1}
\end{equation}
(\cite{M4}). For each $K_\bc$-orbit $S$ on $X$ we defined in \cite{GM1} a subset $C(S)$ of $G_\bc$ by
$$C(S)=\{x\in G_\bc\mid xS\cap S'\mbox{ is non-empty and compact}\}$$
where $S'$ is the $G_\br$-orbit on $X$ given by (1.1).

Akhiezer and Gindikin defined a domain $D/K_\bc$ in $G_\bc/K_\bc$ in \cite{AG} as follows. Let $\frak{g}_\br=\frak{k}\oplus\frak{m}$ denote the Cartan decomposition of $\frak{g}_\br={\rm Lie}(G_\br)$ with respect to $K$. Let $\frak{t}$ be a maximal abelian subspace in $i\frak{m}$. Put
$$\frak{t}^+=\{Y\in\frak{t}\mid |\alpha(Y)|<{\pi\over 2} \mbox{ for all }\alpha\in\Sigma\}$$
where $\Sigma$ is the restricted root system of ${\frak g}_\bc$
with respect to $\frak{t}$. Then $D$ is defined by
$$D=G_\br(\exp\frak{t}^+)K_\bc.$$

We conjectured in \cite{GM1} the following.

\begin{conjecture} \ {\rm (Conjecture 1.6 in \cite{GM1})} \ Suppose that $X=G_\bc/P$ is not $K_\bc$-homogeneous. Then we will have $C(S)_0=D$ for all $K_\bc$-orbits $S$ on $X$ of nonholomorphic type. Here $C(S)_0$ is the connected component of $C(S)$ containing the identity.
\end{conjecture}

\begin{remark} \ (i) \ When $G_\br$ is of Hermitian type, there exist two special closed $K_\bc$-orbits $S_1=K_\bc B/B=Q/B$ and $S_2=K_\bc w_0B/B=w_0Q/B$ on the full flag manifold $G_\bc/B$ where $Q=K_\bc B$ is the usual maximal parabolic subgroup of $G_\bc$ defined by a nontrivial central element in $i\frak{k}$ and $w_0$ is the longest element in the Weyl group. For each parabolic subgroup $P$ containing $B$, two closed $K_\bc$-orbits $S_1P$ and $S_2P$ on $G_\bc/P$ are called of holomorphic type and all the other $K_\bc$-orbits are called of nonholomorphic type. When $G_\br$ is of non-Hermitian type, we define that all the $K_\bc$-orbits are of nonholomorphic type.

(ii) \ If $X=G_\bc/P$ is $K_\bc$-homogeneous, then we have $S=S'=X$ and therefore $C(S)=G_\bc$. So we must assume that $X$ is not $K_\bc$-homogeneous. When $G_\bc$ is simple, it is shown in \cite{O} Theorem 6.1 that there are only two types of $K_\bc$-homogeneous flag manifolds $X$ as follows. (Note that the $K_\bc$-orbit structure on $X$ depends only on the Lie algebras if $K_\bc$ is connected.)

(1) \ $G_\bc=SL(2n,\bc),\ K_\bc=Sp(n,\bc),\ X=P^{2n-1}(\bc).$

(2) \ $G_\bc=SO(2n,\bc),\ K_\bc=SO(2n-1,\bc),\ X=SO(2n)/U(n).$
\end{remark}

Let $S_\op$ denote the unique open $K_\bc$-$B$ double coset in $G_\bc$. Then $S'_\op$ is closed in $G_\bc$ and therefore we can write
$$C(S_\op)=\{x\in G_\bc\mid xS_\op\supset S'_\op\}.$$
The domain $C(S_\op)_0$ is often called the ``Iwasawa domain''.

It is proved by Barchini (\cite{B}) that
$C(S_\op)_0\subset D$.
On the other hand, Huckleberry (\cite{H}) proved the opposite inclusion
\begin{equation}
D\subset C(S_\op)_0. \tag{1.2}
\end{equation}
(Recently \cite{Ms} gave a proof of (1.2) without complex analysis.) So we have the equality
\begin{equation}
C(S_\op)_0=D. \tag{1.3}
\end{equation}

It is proved in Proposition 8.1 and Proposition 8.3 in \cite{GM1} that
$C(S_\op)_0\subset C(S)_0$
for all $K_\bc$-orbits $S$ on all flag manifolds $X=G_\bc/P$. So we have only to prove the inclusion
$C(S)_0\subset D$
in Conjecture 1.1.

The first aim of this paper is the following generalization of Barchini's theorem which solves Conjecture 1.1 for the open $K_\bc$-orbit on an arbitrary flag manifold $G_\bc/P$.

\begin{theorem} \ Suppose that $G_\br$ is simple. Then there exists a $K_\bc$-$B$ invariant subset $\widetilde{S}$ in $G_\bc$ satisfying the following three conditions.

{\rm (i)} \ $\widetilde{S}$ consists of single $K_\bc$-$B$ double coset when $G_\br$ is of non-Hermitian type and consists of two $K_\bc$-$B$ double cosets when $G_\br$ is of Hermitian type.

{\rm (ii)} \ $x\widetilde{S}^{cl}\cap S'_\op\ne \phi$ for all elements $x$ in the boundary of $D$.

{\rm (iii)} \ Let $P$ be a parabolic subgroup of $G_\bc$ containing $B$. If $G_\bc/P$ is not $K_\bc$-homogeneous, then
$$\widetilde{S}\cap S_\op P=\phi.$$
\end{theorem}

\begin{corollary} \ Let $P$ be a parabolic subgroup of $G_\bc$ containing $B$. If $G_\bc/P$ is not $K_\bc$-homogeneous, then $C(S_\op P)_0=D$.
\end{corollary}

\begin{proof} Let $x$ be an element in the boundary of $D$. Then it follows from Theorem 1.3 that $x\widetilde{S}^{cl}\cap S'_\op\ne \phi$. If $x\in C(S_\op P)_0$, then we have
$$xS_\op P\supset S'_\op P.$$
Hence we have $x\widetilde{S}^{cl}\cap xS_\op P\ne \phi$. Since $S_\op P$ is open in $G_\bc$, this implies that
$\widetilde{S}\cap S_\op P\ne\phi$
a contradiction to the condition (iii) in Theorem 1.3. Thus we have $C(S_\op P)_0\subset D$.
\end{proof}

\begin{example} \ Let $G_\bc=SL(3,\bc),\ G_\br=SU(2,1)$ and
$$K_\bc=\left\{\bp * & * & 0 \\ * & * & 0 \\ 0 & 0 & * \ep \in G_\bc \right\}.$$
Then the complex symmetric space $G_\bc/K_\bc$ is identified with the space consisting of the pairs $(V_+,V_-)$ of two-dimensional subspaces $V_+$ and one-dimensional subspaces $V_-$ of $\bc^3$ such that $V_+\cap V_-=\{0\}$ by the identification
$$gK_\bc\mapsto (V_+,V_-)=(gV_+^0, gV_-^0).$$
Here $V_+^0=\bc e_1\oplus \bc e_2$ and $V_-^0=\bc e_3$ with the canonical basis $e_1,e_2,e_3$ of $\bc^3$. Using the Hermitian form
$Q(z)=|z_1|^2+|z_2|^2-|z_3|^2$
defining $SU(2,1)$, we can decompose $\bc^3$ as
$$\bc^3=C_0\sqcup C_+\sqcup C_-$$
where $C_0=\{z\in\bc^3\mid Q(z)=0\},\ C_+=\{z\in\bc^3\mid Q(z)>0\}$ and $C_-=\{z\in\bc^3\mid Q(z)<0\}$. Then the Akhiezer-Gindikin domain $D/K_\bc$ is described as
$$D/K_\bc =\{(V_+,V_-)\in G_\bc/K_\bc\mid V_+-\{0\}\subset C_+ \mbox{ and }V_--\{0\}\subset C_-\}$$
by \cite{GM1} Proposition 2.2. Hence the boundary of $D/K_\bc$ consists of the three $G_\br$-orbits
\begin{align*}
D_1 & =\{(V_+,V_-)\in G_\bc/K_\bc\mid V_+\mbox{ is tangent to }C_0 \mbox{ and }V_--\{0\}\subset C_-\}, \\
D_2 & =\{(V_+,V_-)\in G_\bc/K_\bc\mid V_+-\{0\}\subset C_+\mbox{ and }V_-\in C_0\}, \\
D_3 & =\{(V_+,V_-)\in G_\bc/K_\bc\mid V_+\mbox{ is tangent to }C_0 \mbox{ and }V_-\in C_0\}.
\end{align*}
%(We see that the boundary of $D/K_\bc$ consists of a finite number of $G_\br$-orbits if and only if the real rank of $G_\br$ is one.)

Let $B$ denote the standard Borel subgroup of $G_\bc$ consisting of upper triangular matrices contained in $G_\bc$. Then the full flag manifold $X=G_\bc/B$ consists of flags $(\ell,p)$ where $\ell$ are one-dimensional subspaces of $\bc^3$ and $p$ are two-dimensional subspaces of $\bc^3$ containing $\ell$. Note that $B$ is the isotropy subgroup of the flag $(\bc e_1,\bc e_1\oplus\bc e_2)$. We see that $X$ is decomposed into the six $K_\bc$-orbits
\begin{align*}
S_1 & =\{(\ell,p)\in X\mid \ell=V_-^0\}, \\
S_2 & =\{(\ell,p)\in X\mid p=V_+^0\}, \\
S_3 & =\{(\ell,p)\in X\mid \ell\subset V_+^0\mbox{ and }p\supset V_-^0\}, \\
S_4 & =\{(\ell,p)\in X\mid p\supset V_-^0\}-(S_1\sqcup S_3), \\
S_5 & =\{(\ell,p)\in X\mid \ell\subset V_+^0\}-(S_2\sqcup S_3), \\
S_\op & =X-(S_1\sqcup S_2\sqcup S_3\sqcup S_4\sqcup S_5).
\end{align*}
On the other hand, the corresponding $G_\br$-orbits are
\begin{align*}
S'_1 & =\{(\ell,p)\in X\mid \ell-\{0\}\subset C_-\}, \\
S'_2 & =\{(\ell,p)\in X\mid p-\{0\}\subset C_+\}, \\
S'_3 & =\{(\ell,p)\in X\mid \ell-\{0\}\subset C_+ \mbox{ and }p\cap C_-\ne\phi\}, \\
S'_4 & =\{(\ell,p)\in X\mid \ell\subset C_0\mbox{ and }p\mbox{ is not tangent to }C_0\}, \\
S'_5 & =\{(\ell,p)\in X\mid p\mbox{ is tangent to }C_0\mbox{ and } \ell\not\subset C_0\}, \\
S'_\op & =\{(\ell,p)\in X\mid \ell\subset C_0\mbox{ and }p\mbox{ is tangent to }C_0\},
\end{align*}
respectively.

If $xK_\bc=(V_+,V_-)\in D_1\sqcup D_3$, then $V_+$ is tangent to $C_0$. Hence the flag
$(V_+\cap C_0,V_+)$
is contained in $xS_2\cap S'_\op$. On the other hand, if $xK_\bc=(V_+,V_-)\in D_2\sqcup D_3$, then $V_-\subset C_0$. Hence the flag
$(V_-,p)$
($p$ is tangent to $C_0$) is contained in $xS_1\cap S'_\op$. Thus we have
$$x(S_1\sqcup S_2)\cap S'_\op\ne\phi$$
for all elements $x$ in the boundary of $D$.

There are two nontrivial parabolic subgroups
$$P_1=\{g\in G_\bc\mid gV_+^0=V_+^0\}\mand P_2=\{g\in G_\bc\mid g\bc e_1=\bc e_1\}$$
in $G_\bc$ containing $B$. We see that $S_\op P_1=S_5\sqcup S_\op$ and that $S_\op P_2=S_4\sqcup S_\op$. So we have
$$(S_1\sqcup S_2)\cap S_\op P=\phi$$
for all parabolic subgroups $P$ of $G_\bc$ such that $B\subset P\ne G_\bc$. Thus we have verified Theorem 1.3 for $\widetilde{S}=S_1\sqcup S_2$ in this case.
\end{example}

\begin{remark} \ (i) \ The following statement is false:
$$x\in \partial (C(S)_0)\Longrightarrow xS^{cl}\cap\partial S'\ne\phi$$
for non-open $S$.
In fact there is a counter example when $G_\br=SU(2,1)$ as follows. 
In the above example, let $xK_\bc=(V_+,V_-)$ be a point in $D_1$. Then $V_+$ is tangent to $C_0$ and $V_--\{0\}\subset C_-$. Consider the $K_\bc$-orbit $S_4$ on $G_\bc/B$. Then the intersection $xS_4\cap S'_4$ consists of the flags $(\ell,p)$ such that
$$\ell\subset C_0,\ \ell\notin V_+\mbox{ and } p\supset V_-.$$
Hence $xS_4\cap S'_4$ is not closed in $G_\bc$ and therefore $x\in\partial (C(S_4)_0)$. On the other hand, since $xS_4^{cl}$ consists of flags $(\ell,p)$ satisfying $p\supset V_-$, it does not intersect $\partial S'_4 =S'_\op$.

(ii) \ By the above argument we see that
\begin{equation}
x\in D_1\sqcup D_3\Longrightarrow (xS_4\cap S'_4)^{cl}\supset xS_3\cap S'_4. \tag{1.4}
\end{equation}
On the other hand suppose that $x\in D_2\sqcup D_3$. Then $V_-\subset C_0$. Let $(\ell_0,p_0)$ denote the unique flag in $xS_1\cap S'_\op$ given by
$$\ell_0=V_-\mand p_0\mbox{ is tangent to }C_0.$$
Then we can take a sequence $\{(\ell_n,p_n)\}$ of flags in $xS_4\cap S'_4$ converging to $(\ell_0,p_0)$. Hence
\begin{equation}
x\in D_2\sqcup D_3\Longrightarrow (xS_4\cap S'_4)^{cl}\supset xS_1\cap S'_\op. \tag{1.5}
\end{equation}
By (1.4) and (1.5), we have
$$C(S_4)_0\subset D.$$
Since the opposite inclusion is known by (1.3) and \cite{GM1} Proposition 8.3, we have
$$C(S_4)_0=D.$$
On the other hand, we have $C(S_4P_1)_0=C(S_1)$ since $S_4P_1=S_1P_1$ is of holomorphic type (\cite{WZ1}). Hence $C(S_4)_0$ is strictly included in $C(S_4P_1)_0$ though the contradicting converse inclusion and the equality are asserted in \cite{HN} (ver.2) Proposition 6 and Proposition 10, respectively. (Note that the projection $S'_4/B\to S'_4P_1/P_1$ is proper and that $S_4^{cl}$ is left $P_1$-invariant.)
\end{remark}

Next suppose that $S$ is closed. Then $S'$ is open (\cite{M2}) and so the condition
$$xS\cap S'\mbox{ is non-empty and compact (in $G_\bc/P$)}$$
implies
$$xS\subset S'.$$
Hence the set $C(S)_0$ is the cycle space for $S'$ introduced by Wells and Wolf (\cite{WW}).
Let $B$ be a Borel subgroup of $G_\bc$ contained in $P$. Let $\{S_j \mid j\in J\}$ denote the set of $K_\bc$-$B$ double cosets in $G_\bc$ of codimension one and $T_j=S_j^{cl}$ the closure of $S_j$. Then we defined in \cite{GM2} a subset $J'=J(S)$ of $J$ for $S$ by
$$J'=\{j\in J\mid S(BwB)^{cl}=T_j\mbox{ for some }w\in W\}$$
and proved the equality
\begin{equation}
C(S)_0=\Omega(J') \tag{1.6}
\end{equation}
where
$$\Omega(J')=\{x\in G_\bc\mid xT_j\cap S'_\op=\phi\mbox{ for all }j\in J'\}_0.$$
(See \cite{GM2} Remark 4 for the related paper \cite{HW2}.) It is also shown in \cite{GM2} that
$\bigcup_{j\in J} T_j$
is the complement of $S_\op$ in $G_\bc$. So we have the equalities
\begin{equation}
D=C(S_\op)_0=\Omega(J) \tag{1.7}
\end{equation}
by (1.3).

The second aim of this paper is to prove the following theorem.

\begin{theorem} \ Suppose that $G_\br$ is of non-Hermitian type. Then the $K_\bc$-$B$ double coset $\widetilde{S}$ given in Theorem 1.3 is contained in
$\bigcap_{j\in J} T_j$.
\end{theorem}

As a direct consequence of this theorem, we have:

\begin{corollary} \ If $G_\br$ is of non-Hermitian type, then $\Omega(J')=D$ for all nonempty subsets $J'$ of $J$.
\end{corollary}

\begin{proof} Since $D=\Omega(J)\subset\Omega(J')$ by (1.7), we have only to show the inclusion $\Omega(J')\subset D$. Let $x$ be an arbitrary element on the boudary of $D$. Then
$$xT_j\cap S'_\op\supset x\widetilde{S}^{cl}\cap S'_\op\ne\phi\quad\mbox{for all }j\in J$$
by Theorem 1.3 and Theorem 1.7. So we have
$x\notin \Omega(J')$
for any $J'\ne\phi$. Hence $\Omega(J')\subset D$. \end{proof}

\bigskip
By (1.6) we have:

\begin{corollary} \ If $G_\br$ is of non-Hermitian type and $S$ is a closed $K_\bc$-orbit on an arbitrary flag manifold $X=G_\bc/P$ such that $S\ne X$, then $C(S)_0=D$.
\end{corollary}

\bigskip
Here we review the history of the inclusion
\begin{equation}
C(S)_0\subset D \tag{1.8}
\end{equation}
for closed $K_\bc$-orbits $S$ on $X$. (As we explained, the opposite inclusion is already established.)

It was conjectured in \cite{G} (Problem 6) that the Akhiezer-Gindikin domain $D$ would be the universal domain for all Stein extensions of Riemannian symmetric spaces $G_\br/K$ and that it would coincide with the cycle spaces $C(S)_0$ and the Iwasawa domain $C(S_{\rm op})_0$. 

When $G_\br$ is of Hermitian type, (1.8) was proved in \cite{WZ1} and \cite{WZ2} for closed $K_\bc$-orbits of nonholomorphic type (c.f. Remark 5 in \cite{GM2}). Remark that $C(S)_0$ is bigger than $D$ if $S$ is of holomorphic type (c.f. \cite{WZ1} and \cite{GM1}).

Here we should note the following generality: Let $\pi:X'\to X$ be a $G_\bc$-equivariant surjection between flag manifolds and let $S$ be a closed $K_\bc$-orbit on $X'$. Then
$$xS\subset S'\Longrightarrow x\pi(S)\subset \pi(S').$$
Hence
$$x\in C(S)\Longrightarrow x\in C(\pi(S))$$
and so we have $C(S)\subset C(\pi(S))$. This implies that we have only to prove the inclusion (1.8) for minimal flag manifolds with nontrivial $K_\bc$-orbit structure.

In \cite{GM1}, (1.8) was proved for typical minimal flag manifolds $X$ with nontrivial $K_\bc$-orbit structure for all non-Hermitian classical $G_\br$ by case-by-case checkings. (Note that they include the complex cases that $G_\br=SL(n,\bc),\ SO(n,\bc)$ or $Sp(n,\bc)$.) We proposed our Conjecture 1.1 (Conjecture 1.6 in \cite{GM1}) by these many explicit computations.

Recently, Fels and Huckleberry (\cite{FH}) gave a general proof of (1.8) for closed $S$ for all non-Hermitian cases by using a complex analytic notion ``Kobayashi hyperbolicity''.
But we need no complex analysis (except Lemma 2.1) in our proof of Theorem 1.7, Corollary 1.8 and Corollary 1.9 in this paper.  

The rest of this paper is consructed as follows. In Section 2 and Section 3, we consider real Lie groups $G$ and associated pairs of symmetric subgroups $H$ and $H'$ of $G$. If $G=G_\bc$ and $H=K_\bc$, then $H'=G_\br$. In this general setting we defined a generalization of the Akhiezer-Gindikin domain $D$ in \cite{Ms}. In Section 2 we define generic elements in the boundary of $D$ generalizing the results in \cite{FH} Section 4. (Remark: In Section 3 of \cite{FH}, they studied the Jordan decomposition and elliptic elements of the double coset decomposition $G_\br \backslash G_\bc/K_\bc$. But these results were already given in \cite{M5} in a more general setting as we explain in Section 2.)

In Section 3 we construct a parabolic subgroup $P_Z$ such that $H'P_Z$ and $Ha_\alpha^{-1}P_Z$ are closed in $G$. Here $a_\alpha$ is an element of $T$ contained in the boundary of $\exp\frak{t}^+$. Then we prove the key theorem (Theorem 3.2) which asserts that
$$xHa_\alpha^{-1}P_Z\cap H'P_Z\ne\phi$$
for all elements $x$ in the boundary of $D$ under some conditions. (If $H'$ is a group of non-Hermitian type and $G$ is the complexification of $H'$, then the conditions are satisfied. If $H'$ is of Hermitian type, then we must also consider $P_{-Z}$.)

In sections 4 through 7, we assume that $G$ is a complex semisimple Lie group and $H'$ is a connected real form of $G$. Hence $H$ is the complexification of a maximal compact subgroup of $H'$. (If we use the notations in Section 1, then we rewrite as $G=G_\bc,\ H=K_\bc$ and $H'=G_\br$.) In Section 4 and 5 we prove Theorem 1.3. 
In Section 6 we prove Theorem 1.7. Section 7 is an appendix for the orbit structure on the full flag manifolds.

\begin{remark} Recently Conjecture 1.1 is studied for non-closed and non-open orbits in \cite{Me2}. In this paper our conjecture is solved for non-Hermitian cases. So the remaining problem in our conjecture is only for non-closed and non-open orbits when $G_\br$ is of Hermitian type (c.f. Remark 1.6).
\end{remark}

{\bf Acknowledgement}: The author would like to express his heartily thanks to S. Gindikin for his advice and encouragement.

\section{Generic elements in the boundary of Akhiezer-Gindikin domain}

First we prepare the following lemma on the ``continuity'' of the eigenvalues. For a complex $m\times m$ matrix $A$, define a norm $N(A)$ of $A$ by
$$N(A)=\max\{|a_{ij}|\mid i,j=1,\ldots,m\}.$$
For $\delta>0$ define a compact neighborhood $U_\delta(A)$ of $A$ by
$$U_\delta(A)=\{B\mid N(B-A)\le\delta\}.$$
For a matrix $A$ let $f_A(z)=\det(zI-A)$ denote the eigenpolynomial of $A$.

\begin{lemma} \ {\rm (i)} \ Let $A$ be a complex square matrix with an eigenvalue $\lambda$. Then for any $\varepsilon>0$ there exists a $\delta>0$ such that
$$B\in U_\delta(A) \Longrightarrow \mbox{there exists an eigenvalue $\mu$ of $B$ such that }|\mu-\lambda|<\varepsilon.$$

{\rm (ii)} \ Let $A=\{a_{ij}\}$ be a complex $m\times m$ matrix. Then
$$|\lambda|\le mN(A)$$
for all the eigenvalues $\lambda$ of $A$.
\end{lemma}

\noindent {\it Proof}. (i) \ Take an $\eta$ so that $0<\eta \le \varepsilon$ and that $f_A(z)\ne 0$ for all $z$ on the circle $C_\eta: |z-\lambda|=\eta$. Then there exists a $\delta>0$ such that
$$\ell=\min\{|f_B(z)|\mid B\in U_\delta(A), \ z\in C_\eta\}>0.$$
Put $D_\eta=\{z\in\bc\mid |z-\lambda|\le\eta\}$.

Suppose that there exists a $B\in U_\delta(A)$ such that
$$f_B(z)\ne 0\mbox{ for all }z\in D_\eta.$$
Then we will get a contradiction as follows. Put $B(t)=A+t(B-A)$ for $0\le t\le 1$ and define $t_0$ by
$$t_0=\inf\{t\in [0,1]\mid f_{B(t)}(z)\ne 0\mbox{ for all }z\in D_\eta\}.$$
Then there is a decreasing sequence $\{t_k\}$ in $[t_0,1]$ such that $\lim_{k\to\infty} t_k=t_0$ and that $g_k(z)=1/f_{B(t_k)}(z)$ are holomorphic functions on $D_\eta$. We have
$$|g_k(z)|\le {1\over \ell}$$
for $z\in D_\eta$ by the maximum principle. So we have $|f_{B(t_k)}(z)|\ge \ell$ for $z\in D_\eta$ and therefore $|f_{B(t_0)}(z)|\ge \ell$ for $z\in D_\eta$. Moreover if $t_0>0$, then there exists a $t\in (0,t_0)$ such that $f_{B(t)}(z)\ne 0$ for $z\in D_\eta$, contradicting the definition of $t_0$. So we have $t_0=0$ and $|f_A(z)|\ge \ell$ for $z\in D_\eta$. But this contradicts the assumption that $f_A(\lambda)=0$.

(ii) \ If $|z|>mN(A)$, then we can define
$$(zI-A)^{-1}={1\over z}\left(I-{1\over z}A\right)^{-1}={1\over z}\left(I+{1\over z}A+{1\over z^2}A^2+\cdots\right)$$
since the right hand side converges. So the eigenvalues $\lambda$ of $A$ satisfy
$$\mbx6 |\lambda|\le mN(A).\mbx6 \square$$

\bigskip
Let $G$ be a connected real semisimple Lie group. Let $\sigma$ be an involution of $G$ (automorphism of $G$ of order two). Then there exists a Cartan involution $\theta$ of $G$ such that $\sigma\theta=\theta\sigma$. Denote the corresponding involutions of $\frak{g}={\rm Lie}(G)$ by the same letters as usual. Let
$$\frak{g}=\frak{h}\oplus\frak{q} \mand \frak{g}=\frak{k}\oplus\frak{m}$$
denote the $+1$ and $-1$ eigenspace decomposition of $\frak{g}$ for the involutions $\sigma$ and $\theta$, respectively. Put $\tau=\sigma\theta$. Then $\tau$ is also an involution of $G$. Let
$$\frak{g}=\frak{h}'\oplus\frak{q}'$$
be the $+1$ and $-1$ eigenspace decomposition of $\frak{g}$ for $\tau$. Then
$$\frak{h}'=(\frak{k}\cap\frak{h})\oplus(\frak{m}\cap\frak{q}) \mand
\frak{q}'=(\frak{k}\cap\frak{q})\oplus(\frak{m}\cap\frak{h}).$$
Define two symmetric subgroups $H$ and $H'$ as the connected components of $G^\sigma$ and $G^\tau$, respectively, containing the identity. Then $H$ and $H'$ are called ``associated'' (\cite{Be}, \cite{M1}).

It is studied in \cite{M5} the double coset decomposition $H\backslash G/L$ where $H$ and $L$ are arbitrary two symmetric subgroups of $G$. Of course we can apply it to the pair of the symmetric subgroups $(H,H')$.

\begin{remark} \ $K_\bc$ and $G_\br$ are associated in the complex Lie group $G_\bc$. So we can consider the setting in Section 1 as a special one.
\end{remark}

We defined a generalization of the Akhiezer-Gindikin domain in \cite{Ms} as follows. Let $\frak{t}$ be a maximal abelian subspace of $\frak{k}\cap\frak{q}$. Then we can define the root space
$$\frak{g}_\bc(\frak{t},\alpha)=\{X\in\frak{g}_\bc\mid [Y,X]=\alpha(Y)X\mbox{ for all }Y\in\frak{t}\}$$
for any linear form $\alpha:\frak{t}\to i\br$. Here $\frak{g}_\bc=\frak{g}\oplus i\frak{g}$ is the complexification of $\frak{g}$. Put
$$\Sigma=\Sigma(\frak{g}_\bc,\frak{t})=\{\alpha\in i\frak{t}^*-\{0\} \mid \frak{g}_\bc(\frak{t},\alpha)\ne\{0\}\}.$$
Then $\Sigma$ satisfies the axiom of the root system (\cite{R} Theorem 5). Since $\theta(Y)=Y$ for all $Y\in\frak{t}$, we can decompose $\frak{g}_\bc(\frak{t},\alpha)$ into the $+1,-1$-eigenspaces for $\theta$ as
\begin{equation}
\frak{g}_\bc(\frak{t},\alpha)=\frak{k}_\bc(\frak{t},\alpha) \oplus\frak{m}_\bc(\frak{t},\alpha). \tag{2.1}
\end{equation}
Define a subset
\begin{equation}
\Sigma(\frak{m}_\bc,\frak{t})=\{\alpha\in i\frak{t}^*-\{0\} \mid \frak{m}_\bc(\frak{t},\alpha)\ne\{0\}\} \tag{2.2}
\end{equation}
of $\Sigma$ and put
$$\frak{t}^+=\{Y\in{\frak t}\mid |\alpha(Y)|<{\pi\over 2} \mbox{ for all }\alpha\in\Sigma(\frak{m}_\bc,\frak{t})\}.$$
Then we define a generalization of the Akhiezer-Gindikin domain $D$ in $G$ by
$$D=H'T^+H.$$
where $T^+=\exp\frak{t}^+$. (We showed in \cite{Ms} Proposition 1 that $D$ is open in $G$.) The following assertion is already proved in the proof of \cite{Ms} Proposition 2.

\begin{lemma} \ $D\cap T=T^+(T\cap H)$.
\end{lemma}

In \cite{M5} we considered the automorphism
$$f_x=\tau\Ad(x)\sigma\Ad(x)^{-1}$$
of $\frak{g}$ for every element $x$ in $G$. Then we defined the ``Jordan decomposition'' $x=(\exp X_n)x_s$ of the element $x$ so that $f_x=f_{x_s}\Ad(\exp(-2X_n))=\Ad(\exp(-2X_n))f_{x_s}$ is the usual multiplicative Jordan decomposition of the automorphism $f_x$ of $\frak{g}$ and that $X_n$ is a nilpotent element in $\frak{q}'\cap \Ad(x_s)\frak{q}$. It is shown in \cite{M5} Proposition 3 that
\begin{equation}
H'x_sH\subset(H'xH)^{cl}. \tag{2.3}
\end{equation}
Let $\partial D$ denote the boundary of $D$ in $G$ and $\partial T^+$ the boundary of $T^+$ in $T$.

\begin{lemma} \ If $x\in\partial D$, then $H'x_sH=H'aH$ for some $a\in\partial T^+$.
\end{lemma}

\begin{proof} If $y=h'th\in D$ with $h'\in H',\ t\in T^+$ and $h\in H$, then
$$f_y=\tau\Ad(h'th)\sigma\Ad(h'th)^{-1}=\Ad(h')\tau\sigma\Ad(t)^{-2}\Ad(h')^{-1}=\Ad(h')\theta\Ad(t)^{-2}\Ad(h')^{-1}.$$
Hence the absolute values of the eigenvalues of $f_y$ are all equal to one.

Let $x$ be an element of $\partial D$. By the continuity of the eigenvalues shown in Lemma 2.1 (i), the absolute values of the eigenvalues of $f_x$ are all one. This holds also for $f_{x_s}$. So if we decompose $x_s=(\exp X_p)x_k$ (polar decomposition) as in \cite{M5} Section 4.2, then we have $X_p=0$ and $x_k=x_s$. Let us call such an element $x_s$ ``elliptic'' following the terminology in \cite{FH}. We can show that $y\in G$ is elliptic if and only if $y$ is contained in $H'TH$ by the arguments in the proof of \cite{M5} Theorem 2. Write $x_s=h'ah$ with some $h'\in H',\ a\in T$ and $h\in H$. By (2.3) $x_s$ is contained in $D^{cl}$. Hence $a$ is also contained in $D^{cl}$. We have only to show $a\in (D\cap T)^{cl}$ in view of Lemma 2.3.

Since $f_a=\tau\Ad(a)\sigma\Ad(a)^{-1}$ is semisimple, we can decompose $\frak{g}$ as
\begin{equation}
\frak{g}=(\frak{h}'+\Ad(a)\frak{h})\oplus (\frak{q}'\cap\Ad(a)\frak{q}) \tag{2.4}
\end{equation}
by \cite{M5} Lemma 1. Let $B(\ ,\ )$ denote the Killing form on $\frak{g}$ and let $B_\theta(\ ,\ )$ denote the positive definite bilinear form on $\frak{g}$ defined by $B_\theta(X,Y)=-B(X,\theta Y)$. Since $a\in D^{cl}$, we can take by (2.4) an element $Y$ of $\frak{q}'\cap\Ad(a)\frak{q}$ such that $y=(\exp Y)a\in D$ and that $B_\theta(Y,Y)<\varepsilon^2$ for an arbitrary positive real number $\varepsilon$. Since
$$f_y=\tau\Ad(y)\sigma\Ad(y)^{-1}=f_a\Ad(\exp 2Y)^{-1}=\Ad(\exp 2Y)^{-1}f_a$$
is elliptic, it follows that $\Ad(\exp 2Y)$ is elliptic. Taking $\varepsilon$ sufficiently small, we may assume that $Y$ is semisimple and that the eigenvalues of $\ad(Y)$ are pure imaginary. Hence there exists an element $h\in H'\cap\Ad(a)H$ such that $Y'=\Ad(h)Y$ is contained in $\frak{t}$ which is a compact Cartan subset of $\frak{q}'\cap\Ad(a)\frak{q}$ (c.f. \cite{OM} Corollary of Theorem 2). Since $(\exp Y')a\in D\cap T$ and since $B_\theta(Y',Y')\le B_\theta(Y,Y)<\varepsilon^2$ by the following lemma, we have proved $a\in(D\cap T)^{cl}$.
\end{proof}

\begin{lemma} \ If $Y=\Ad(g)Y'$ for some $g\in G$ and $Y'\in\frak{k}$, then $B_\theta(Y',Y')\le B_\theta(Y,Y)$.
\end{lemma}

\noindent {\it Proof}. Since $B_\theta(\Ad(k)Y',\Ad(k)Y')=B_\theta(Y',Y')$ for $k\in K$, we may assume $g=\exp Z$ for some $Z\in\frak{m}$. Write
$Y'=\sum_\lambda Y_\lambda$
with $\lambda$-eigenvectors $Y_\lambda$ for $\ad(Z)$ ($\lambda\in\br$). Since $Y'\in\frak{k}$, we have $\theta Y_\lambda=Y_{-\lambda}$. Hence
\begin{align*}
\mbx{1.5} B_\theta(Y,Y) & =B_\theta(\Ad(\exp Z)Y',\Ad(\exp Z)Y') 
=\sum_\lambda e^{2\lambda}B_\theta(Y_\lambda,Y_\lambda) \\
& =B_\theta(Y_0,Y_0)+\sum_{\lambda>0} (e^{2\lambda}+e^{-2\lambda})B_\theta(Y_\lambda,Y_\lambda) \\
& \ge B_\theta(Y_0,Y_0)+\sum_{\lambda>0} 2B_\theta(Y_\lambda,Y_\lambda) 
=B_\theta(Y',Y') \mbx3 \square
\end{align*}

\begin{lemma} \ Every $x\in\partial D$ is contained in the boundary of the complement $G-D^{cl}$ of $D^{cl}$ in $G$.
\end{lemma}

\begin{proof} First assume that $x$ is semisimple. Then we may assume $x=a=\exp Y$ with $Y\in \partial \frak{t}^+$ by Lemma 2.4. There exists an $\alpha\in\Sigma$ such that $\alpha(Y)=\pi i/2$ and that $\frak{m}_\bc (\frak{t},\alpha)\ne\{0\}$ by the definition of $\frak{t}^+$. Since $X\mapsto \tau\overline{X}$ defines a conjugation of $\frak{m}_\bc (\frak{t},\alpha)$, we can take a nonzero element $X$ of $\frak{m}_\bc (\frak{t},\alpha)$ such that $\tau\overline{X}=-X$. Then $Z=X+\overline{X}$ is a nonzero element of $\frak{g}^{f_a}\cap\frak{m}\cap\frak{q}'$. Since
\begin{align*}
f_{(\exp tZ)a} & =\tau\Ad(\exp tZ)\Ad(a)\sigma\Ad(a)^{-1}\Ad(\exp tZ)^{-1} \\
& =\Ad(\exp tZ)^{-1}\tau\Ad(a)\sigma\Ad(a)^{-1}\Ad(\exp tZ)^{-1} \\
& =\Ad(\exp tZ)^{-1}f_a\Ad(\exp tZ)^{-1} 
=f_a\Ad(\exp tZ)^{-2} 
=\Ad(\exp tZ)^{-2}f_a
\end{align*}
and since $\Ad(\exp tZ)^{-2}$ has a nontrivial real eigenvalue for $t\in\br^\times$, it follows that $(\exp tZ)a$ is a non-elliptic semisimple element and hence $(\exp tZ)a\notin D^{cl}$ for $t\in\br^\times$. Thus $a\in\partial (G-D^{cl})$.

Next assume that $x$ is not semisimple. Then we may assume $x=(\exp X_n)a$ with an $a\in \partial T^+$ and a nilpotent element $X_n\ne 0$ in $\frak{q}'\cap\Ad(a)\frak{q}$ by Lemma 2.4. By a generalization of the Jacobson-Morozov theorem, there exist a $Y\in\frak{h}'\cap\Ad(a)\frak{h}$ and an $X'_n\in \frak{q}'\cap\Ad(a)\frak{q}$ such that
$$[Y,X_n]=2X_n,\quad [Y,X'_n]=-2X'_n\mand [X_n,X'_n]=Y$$
(\cite{KR} Proposition 4, c.f. also \cite{Se} Section 1). For $t\in\br^\times$ we have
$$[X_n+t^2X'_n,Y-{1\over t}X_n+tX'_n]=2t(Y-{1\over t}X_n+tX'_n).$$
Hence $\Ad(\exp(X_n+t^2X'_n))$ is semisimple and it has an eigenvalue $e^{2t}$. From the same argument as in the first case it follows that $(\exp(X_n+t^2X'_n))a$ is a non-elliptic semisimple element for $t\in\br^\times$. Hence $x=(\exp X_n)a\in\partial (G-D^{cl})$.
\end{proof}

Let $\Sigma(\frak{m}_\bc,\frak{t})$ be as in (2.2) and similarly define another subset $\Sigma(\frak{k}_\bc,\frak{t})$ of $\Sigma$ by
$$\Sigma(\frak{k}_\bc,\frak{t})=\{\alpha\in i\frak{t}^*-\{0\} \mid \frak{k}_\bc(\frak{t},\alpha)\ne\{0\}\}.$$
For $\alpha\in\Sigma$ and $k\in\br$ let $p_{\alpha,k}$ denote the hyperplane in $\frak{t}$ defined by
$$p_{\alpha,k}=\{Y\in\frak{t}\mid \alpha(Y)=k\pi i\}.$$
Define families $\mathcal{H}_+$ and $\mathcal{H}_-$ of hyperplanes by
$$\mathcal{H}_+=\{p_{\alpha,k}\mid \alpha\in\Sigma(\frak{k}_\bc,\frak{t}),\ k\in\bz\}
\mand
\mathcal{H}_-=\{p_{\alpha,k}\mid \alpha\in\Sigma(\frak{m}_\bc,\frak{t}),\ k\in {1\over 2}+\bz\},$$
respectively. Put $\mathcal{H}=\mathcal{H}_+\sqcup \mathcal{H}_-$. Then the family $\mathcal{H}$ of hyperplanes defines a family $\mathcal{E}$ of Euclidean cells consisting of cells $\Delta$ which are maximal connected subsets of $\frak{t}$ satisfying the condition
$$\Delta\subset p\mbox{ or }\Delta\cap p=\phi\mbox{ for all }p\in\mathcal{H}.$$
We have the cellular decomposition
$\frak{t}=\bigsqcup_{\Delta\in\mathcal{E}} \Delta$.
Since $\frak{t}^+$ is bounded (\cite{Ms} Lemma 1), there exists a finite subset $\mathcal{B}$ of $\mathcal{E}$ such that
$\partial \frak{t}^+=\bigsqcup_{\Delta\in\mathcal{B}} \Delta$.
For $\Delta\in\mathcal{E}$ let $\Sigma_\Delta^\pm$ denote the subset of $\Sigma$ given by
$$\Sigma_\Delta^\pm=\{\alpha\in\Sigma\mid \Delta\subset p_{\alpha,k}\mbox{ for some }p_{\alpha,k}\in\mathcal{H}_\pm\}$$
and put $\Sigma_\Delta=\Sigma_\Delta^+\sqcup \Sigma_\Delta^-$.

Let $Y$ be an element of a cell $\Delta$ and put $a=\exp Y$. We can prove the following lemma by the same argument as in \cite{Ms} Lemma 2.

\begin{lemma} \ $\frak{g}_\bc^{f_a}=\frak{z}_{\frak{k}_\bc} (\frak{t})\oplus \bigoplus_{\alpha\in\Sigma_\Delta^+}\frak{k}_\bc(\frak{t},\alpha) \oplus \bigoplus_{\alpha\in\Sigma_\Delta^-}\frak{m}_\bc(\frak{t},\alpha)$.
\end{lemma}

\begin{proof} By (2.1) we have the decomposition
$$\frak{g}_\bc=\bigoplus_{\alpha\in\Sigma\sqcup\{0\}} (\frak{k}_\bc(\frak{t},\alpha)\oplus \frak{m}_\bc(\frak{t},\alpha)).$$
If $X\in\frak{k}_\bc(\frak{t},\alpha)$, then
$f_a(X)=\Ad(a)^{-2}\theta(X)=e^{-2\alpha(Y)}X$. 
If $X\in\frak{m}_\bc(\frak{t},\alpha)$, then
$f_a(X)=\Ad(a)^{-2}\theta(X)=-e^{-2\alpha(Y)}X$. 
So the assertion is clear.
\end{proof}

This lemma implies that
$$\frak{l}_\Delta=\frak{g}^{f_a}=(\frak{h}'\cap\Ad(a)\frak{h})\oplus (\frak{q}'\cap\Ad(a)\frak{q})$$
is independent of the choice of $Y$ in $\Delta$. Let $\frak{z}_\Delta$ denote the center of $\frak{l}_\Delta$ and $\frak{s}_\Delta=[\frak{l}_\Delta,\frak{l}_\Delta]$ the semisimple part of $\frak{l}_\Delta$. Then we have
$\frak{l}_\Delta=\frak{z}_\Delta\oplus\frak{s}_\Delta$
and therefore
\begin{equation}
\frak{q}'\cap\Ad(a)\frak{q}=\frak{l}_\Delta\cap\frak{q}'= (\frak{z}_\Delta\cap\frak{q}') \oplus(\frak{s}_\Delta\cap\frak{q}'). \tag{2.5}
\end{equation}
By Lemma 2.7 we have
$$\frak{z}_\Delta\cap\frak{q}'\subset \frak{z}_\frak{k}(\frak{t})\cap\frak{q}'=\frak{t}.$$
Hence we can write
$$\frak{z}_\Delta\cap\frak{q}'=\{Y\in\frak{t}\mid \alpha(Y)=0\mbox{ for all }\alpha\in\Sigma_\Delta\}.$$
This implies that
\begin{equation}
\frak{z}_\Delta\cap\frak{q}'\mbox{ is the tangent space of }\Delta. \tag{2.6}
\end{equation}

We see that the nilpotent variety in $\frak{s}_\Delta\cap\frak{q}'$ is decomposed into a finite number of $(S_\Delta\cap H')_0$-orbits where $S_\Delta$ denotes the analytic subgroup of $G$ for $\frak{s}_\Delta$. (It is shown in \cite{KR} Theorem 2 that there are a finite number of nilpotent $((S_\Delta\cap H')_0)_\bc$-orbits $\mathcal{N}_\bc$ in $(\frak{s}_\Delta\cap\frak{q}')_\bc$. For each $\mathcal{N}_\bc$ it follows from the Whitney's theorem for real algebraic varieties (\cite{Wh}, \cite{PR} Theorem 3.3) that $\mathcal{N}_\bc\cap(\frak{s}_\Delta\cap\frak{q}')$ consists of a finite number of connected components which are $(S_\Delta\cap H')_0$-orbits.)

For $\Delta\in\mathcal{B}$ let $\mathcal{N}_d$ denote the union of the nilpotent $(S_\Delta\cap H')_0$-orbits of codimension $d$ which are contained in the closure of the set $(\frak{s}_\Delta\cap\frak{q}')_{ell}$ of the elliptic elements in $\frak{s}_\Delta\cap\frak{q}'$. Define a subset
$$M(\Delta,d)=H'(\exp \mathcal{N}_d)(\exp\Delta)H$$
of $G$.

\begin{proposition} \ {\rm (i)} \ $\partial D=\bigsqcup_d \bigcup_{\Delta\in\mathcal{B}} M(\Delta,d)$.

{\rm (ii)} \ If $\mathcal{N}_d\ne\phi$, then $M(\Delta,d)$ is a locally closed $d$-codimensional submanifold of $G$ consisting of a finite number of connected components.

{\rm (iii)} \ $\partial D=(\bigcup_{\Delta\in\mathcal{B}}  M(\Delta,1))^{cl}$.
\end{proposition}

\begin{proof} (i) \ Let $x$ be an element of $\partial D$. By Lemma 2.4 $x$ is contained in
$$H'(\exp X_n)aH$$
for some $a\in\partial T^+$ and a nilpotent element $X_n$ in $\frak{s}_\Delta\cap\frak{q}'$. Here we write $a=\exp Y$ with $Y\in\Delta\subset\partial\frak{t}^+$. We will show that $X_n$ is contained in the closure of $(\frak{s}_\Delta\cap\frak{q}')_{ell}$. Since $S_\Delta\cap H'\subset H'\cap aHa^{-1}$, we can take an $S_\Delta\cap H'$-conjugate of $X_n$ so that $X_n$ is sufficiently close to the origin. If $X_n$ is not on the boundary of $(\frak{s}_\Delta\cap\frak{q}')_{ell}$, then there exists a neighborhood $U$ of $X_n$ in $\frak{s}_\Delta\cap\frak{q}'$ consisting of non-elliptic elements. Take a neighborhood $V$ of $0$ in $\frak{z}_\Delta\cap\frak{q}'$. Then it follows from (2.4) and (2.5) that the set
$$H'(\exp U)(\exp V)aH$$
contains a neighborhood of $y=(\exp X_n)a$ in $G$ consisting of non-elliptic elements. Hence $y$ is not contained in the closure of $D$, contradicting to $x\in\partial D$. Thus we have proved that $X_n$ is on the boundary of $(\frak{s}_\Delta\cap\frak{q}')_{ell}$. Conversely, if $X_n$ is on the boundary of $(\frak{s}_\Delta\cap\frak{q}')_{ell}$, then it is clear that $y=(\exp X_n)a\in\partial D$.

(ii) \ It follows from (2.4), (2.5) and (2.6) that the codimension of $M(\Delta,d)$ in $G$ equals the codimension $d$ of $\mathcal{N}_d$ in $\frak{s}_\Delta\cap\frak{q}'$. 
Considering the left $H'$-action and the right $H$-action, we have only to show that $M(\Delta,d)$ is locally closed at $y=\exp X_n\exp Y_0$ for every $Y_0\in\Delta$ and every $X_n\in \mathcal{N}_d$. Furthermore taking an $S_\Delta\cap H'$-conjugate, we may assume that $X_n$ is sufficiently close to $0$. Let $V$ be a compact neighborhood of $Y_0$ in $\Delta$. Then the nontrivial eigenvalues of $f_{\exp Y}$ for $Y\in V$ are contained in a compact subset $\Lambda$ of $U(1)=\{z\in\bc\mid |z|=1\}$ such that $1\notin\Lambda$. By Lemma 2.1 (ii) we can take a neighborhood $U$ of $0$ in $\frak{s}_\Delta\cap\frak{q}'$ such that the eigenvalues of $\Ad\exp(-2Z)$ are not contained in $\Lambda^{-1}$ for all $Z\in U$. By (2.4) and (2.5) we have only to show that
$$(\exp U\exp V)\cap M(\Delta,d)=\exp(U\cap\mathcal{N}_d)\exp V.$$
Suppose $Z\in U,\ Y\in V$ and $\exp Z\exp Y\in M(\Delta,d)$. Then we have only to show $Z\in\mathcal{N}_d$. 

Let $\frak{g}_\bc^\lambda$ denote the $\lambda$-eigenspace for $f_{\exp Y}$. Then $\frak{g}_\bc$ is decomposed as
$$\frak{g}_\bc=\bigoplus_\lambda \frak{g}_\bc^\lambda$$
since $f_{\exp Y}$ is semisimple. Since $f_{\exp Z\exp Y}=f_{\exp Y}\Ad\exp(-2Z)=\Ad\exp(-2Z)f_{\exp Y}$, we can furthurmore decompose the spaces $\frak{g}_\bc^\lambda$ into the generalized eigenspaces as
$$\frak{g}_\bc^\lambda=\bigoplus_\mu \frak{g}_\bc^{\lambda,\mu}$$
where $\frak{g}_\bc^{\lambda,\mu}=\{X\in\frak{g}_\bc^\lambda\mid (\Ad\exp(-2Z)-\mu I)^kX=0\mbox{ for some }k\}$. Note that every $X\in\frak{g}_\bc^{\lambda,\mu}$ is a generalized eigenvector for $f_{\exp Z\exp Y}$ with the eigenvalue $\lambda\mu$. Since we assume that $\mu\notin\Lambda^{-1}$, we have
$$\lambda\mu=1\Longrightarrow \lambda=\mu=1.$$
This implies that the generalized eigenspace $\frak{g}_\bc^{(1)}$ of $f_{\exp Z\exp Y}$ for the eigenvalue $1$ is contained in $\frak{g}_\bc^1=(\frak{l}_\Delta)_\bc$. On the other hand, the condition $\exp Z\exp Y\in M(\Delta,d)=H'(\exp \mathcal{N}_d)(\exp\Delta)H$ implies that
\begin{equation}
f_{\exp Z\exp Y}\mbox{ is conjugate to }f_{(\exp W)b} \tag{2.7}
\end{equation}
for some $W\in\mathcal{N}_d$ and $b\in\exp\Delta$. So we have
$$\dim\frak{g}_\bc^{(1)}=\dim\frak{g}_\bc^{f_b}=\dim(\frak{l}_\Delta)_\bc.$$
Hence we have
$$\frak{g}_\bc^{(1)}= (\frak{l}_\Delta)_\bc.$$
This implies that the linear map $\Ad\exp(-2Z): (\frak{l}_\Delta)_\bc\to(\frak{l}_\Delta)_\bc$ has the unique eigenvalue $1$. Hence $\Ad\exp(-2Z)$ is unipotent on $(\frak{l}_\Delta)_\bc$ and therefore $Z\in\frak{s}_\Delta\cap\frak{q}'$ is nilpotent.

It follows from (2.7) that
$$\dim_\bc \frak{g}_\bc^{f_{\exp Z\exp Y}}=\dim_\bc \frak{g}_\bc^{f_{(\exp W)b}}.$$
Since $\frak{g}_\bc^{f_{\exp Z\exp Y}}$ is contained in $\frak{g}_\bc^{(1)}=\frak{g}_\bc^1= (\frak{l}_\Delta)_\bc$, we have
$$\frak{g}_\bc^{f_{\exp Z\exp Y}}=(\frak{l}_\Delta)_\bc^{\exp Z}= (\frak{l}_\Delta)_\bc^Z=\{X\in\frak{l}_\bc\mid [X,Z]=0\}.$$
Similary we also have
$$\frak{g}_\bc^{f_{(\exp W)b}}= (\frak{l}_\Delta)_\bc^W=\{X\in\frak{l}_\bc\mid [X,W]=0\}.$$
Hence we have $\dim_\bc (\frak{l}_\Delta)_\bc^Z=\dim_\bc (\frak{l}_\Delta)_\bc^W$ and therefore $\dim_\bc (\frak{s}_\Delta)_\bc^Z=\dim_\bc (\frak{s}_\Delta)_\bc^W$. This implies that
$$\dim_\bc \Ad(S_\Delta)_\bc Z=\dim_\bc \Ad(S_\Delta)_\bc W.$$
By \cite{KR} Proposition 5 we have
$$\dim_\bc \Ad(S_\Delta\cap H')_\bc Z
={1\over 2}\dim_\bc \Ad(S_\Delta)_\bc Z
={1\over 2}\dim_\bc \Ad(S_\Delta)_\bc W
=\dim_\bc \Ad(S_\Delta\cap H')_\bc W.$$
So we have
\begin{equation}
\dim\,\Ad(S_\Delta\cap H') Z=\dim\,\Ad(S_\Delta\cap H') W. \tag{2.8}
\end{equation}
Since $W\in\mathcal{N}_d$, it is contained in the closure of $(\frak{s}_\Delta\cap\frak{q}')_{ell}$. Hence $\exp Z\exp Y\in H'(\exp W)aH$ is contained in $\partial D$. By the argument in the proof of (i), we see that $Z$ is contained in the closure of $(\frak{s}_\Delta\cap\frak{q}')_{ell}$. Combining with (2.8), we have proved $Z\in\mathcal{N}_d$.

(iii) \ Let $x$ be an element of $M(\Delta,d)\subset \partial D$ such that $d\ge 2$. Let $U$ be a neighborhood of $x$ in $G$. Suppose that $U$ does not intersect $\bigcup_{\Delta\in\mathcal{B}} M(\Delta,1)$. Then $U-\partial D$ is connected because $U\cap \partial D$ is a finite union of locally closed submanifolds of codimension greater than two. But $U\cap D$ and $U-D^{cl}$ are both nonempty open sets by Lemma 2.6, a contradiction. Thus $\partial D=(\bigcup_{\Delta\in\mathcal{B}} M(\Delta,1))^{cl}$.
\end{proof}

\begin{remark} \ (i) \ It is known that $d\ge\dim(\frak{s}_\Delta\cap\frak{t})$ (\cite{KR} Proposition 9). Hence
$$M(\Delta,1)\ne\phi$$
only when the rank of $\Sigma_\Delta=\Sigma_\Delta^+\sqcup \Sigma_\Delta^-$ is one.

(ii) \ In general we can prove that ${\rm codim}_G H'xH\ge \dim\frak{t}$ for any $x=(\exp X_n)x_s\in G$. So we may call $x=(\exp X_n)x_s$ regular if ${\rm codim}_G H'xH=\dim\frak{t}$. Regular semisimple elements studied in \cite{M5} are contained in the set of the regular elements. The author is grateful to Michael Otto for suggesting the existence of this notion of regularity.

(iii) \ Suppose that $\frak{m}_\bc(\frak{t},\alpha)\ne 0$ for all the longest roots $\alpha$ in $\Sigma$. Since all the roots in $\Sigma$ are contained in the convex hull spanned by the longest roots, $\frak{t}^+$ is written as
$$\frak{t}^+=\{Y\in\frak{t}\mid |\alpha(Y)|<{\pi\over 2}\mbox{ for all the longest roots }\alpha\in\Sigma\}.$$
So the condition $M(\Delta,1)\ne\phi$ for $\Delta\in\mathcal{B}$ implies that
$$\Sigma_\Delta=\Sigma_\Delta^-=\{\pm\alpha\}$$
for some longest root $\alpha$. Suppose furthermore that all the longest roots in $\Sigma$ are mutually $W_{K\cap H}(\frak{t})$-conjugate. Then we can write
$$\bigcup_{\Delta\in\mathcal{B}} M(\Delta,1)=\bigcup_{\Delta\in\mathcal{B}_\alpha} M(\Delta,1)$$
with a longest root $\alpha$ in $\Sigma$ where
$$\mathcal{B}_\alpha=\{\Delta\in\mathcal{B}\mid \alpha(Y)={\pi\over 2}i\mbox{ for all }Y\in\Delta\}.$$
So we have
$\partial D=\left(\bigcup_{\Delta\in\mathcal{B}_\alpha} M(\Delta,1)\right)^{cl}.$
\end{remark}

\section{Construction of parabolic subgroups}

Let $\alpha$ be a root in $\Sigma(\frak{m}_\bc,\frak{t})$ and suppose that $2\alpha\notin\Sigma$. Let $\frak{s}_\bc$ be the subalgebra of $\frak{g}_\bc$ generated by $\frak{m}_\bc(\frak{t},\alpha) \oplus \frak{m}_\bc(\frak{t},-\alpha)$. Since
$[\frak{m}_\bc(\frak{t},\alpha),\frak{m}_\bc(\frak{t},-\alpha)] \subset \frak{z}_{\frak{k}_\bc}(\frak{t})$
and since
$[\frak{z}_{\frak{k}_\bc}(\frak{t}), \frak{m}_\bc(\frak{t},\pm\alpha)] \subset \frak{m}_\bc(\frak{t},\pm\alpha),$
we have
\begin{equation}
\frak{s}_\bc\subset \frak{z}_{\frak{k}_\bc}(\frak{t})\oplus \frak{m}_\bc(\frak{t},\alpha) \oplus \frak{m}_\bc(\frak{t},-\alpha). \tag{3.1}
\end{equation}
Clearly $\frak{s}=\frak{s}_\bc\cap\frak{g}$ is a real form of $\frak{s}_\bc$.

\begin{proposition} \ $(\frak{s},\ \frak{s}^\tau)\cong (\frak{so}(2,n),\ \frak{so}(1,n))$ where $n=\dim_\bc\frak{m}_\bc (\frak{t},\alpha)$.
\end{proposition}

\begin{proof} Let $X\mapsto \overline{X}$ denote the conjugation with respect to the real form $\frak{g}$. Since $X\mapsto \overline{\tau(X)}$ is a conjugation of $\frak{m}_\bc(\frak{t},\alpha)$, we can take a nonzero element $X\in\frak{m}_\bc(\frak{t},\alpha)$ such that $\tau(X)=\overline{X}$. Put $Y_1=[X,\tau(X)]=[X,\overline{X}]$. Then we have
$$Y_1\in \frak{z}_{\frak{g}_\bc}(\frak{t})\cap \frak{k}_\bc\cap \frak{q}'_\bc\cap i\frak{g}=\frak{z}_{\frak{g}_\bc}(\frak{t})\cap i(\frak{k}\cap \frak{q})=i\frak{t}.$$
On the other hand, we have
$$B(Y,Y_1)=B(Y,[X,\overline{X}])= B([Y,X],\overline{X})= \alpha(Y)B(X,\overline{X})$$
for $Y\in\frak{t}_\bc$. Since the Hermitian form $B(X,\overline{X})$ is positive definite on $\frak{m}_\bc$, we see that $Y_1\ne 0$. Taking $Y=Y_1$, we have
\begin{equation}
\alpha(Y_1)>0. \tag{3.2}
\end{equation}
Since $Y_1\in\frak{s}_\bc$ and since $[Y_1,\ \frak{m}_\bc(\frak{t},\pm\alpha)]=\frak{m}_\bc(\frak{t},\pm\alpha)$ by (3.2), the spaces $\frak{m}_\bc(\frak{t},\pm\alpha)$ are contained in the derived ideal $[\frak{s}_\bc,\frak{s}_\bc]$ of $\frak{s}_\bc$. So we have $\frak{s}_\bc=[\frak{s}_\bc,\frak{s}_\bc]$ and therefore $\frak{s}_\bc$ is semisimple. If $Z\in\frak{s}_\bc\cap\frak{t}_\bc$ satisfies $\alpha(Z)=0$, then $Z$ is contained in the center of $\frak{s}_\bc$ and hence $Z=0$. So we have proved
\begin{equation}
\frak{s}_\bc\cap\frak{t}_\bc=\bc Y_1. \tag{3.3}
\end{equation}

We will show that $\bc Y_1$ is maximal abelian in $\frak{s}_\bc\cap \frak{q}'_\bc$. Let $X$ be an element in $\frak{s}_\bc\cap \frak{q}'_\bc$ such that $[Y_1,X]=0$. If $Y\in\frak{t}$ satisfies $\alpha(Y)=0$, then $[Y,X]=0$ because $\frak{s}_\bc$ is generated by $\frak{m}_\bc(\frak{t},\alpha) \oplus \frak{m}_\bc(\frak{t},-\alpha)$. Hence $X$ is contained in the centralizer of $\frak{t}$. By (3.1) $X$ is contained in $\frak{z}_{\frak{k}_\bc}(\frak{t})$. Hence $X\in\frak{k}_\bc\cap\frak{q}'_\bc=\frak{k}_\bc\cap \frak{q}_\bc$.  
Since $\frak{t}$ is maximal abelian in $\frak{k}\cap \frak{q}$, we have $X\in\frak{t}_\bc$. It follows from (3.3) that $X\in\bc Y_1$.
 
Thus the symmetric pair $(\frak{s}_\bc,\frak{s}_\bc\cap\frak{h}'_\bc)$ is rank one and we have the restricted root space decomposition
$$\frak{s}_\bc=\frak{z}_{\frak{s}_\bc}(\frak{t})\oplus \frak{m}_\bc(\frak{t},\alpha)\oplus \frak{m}_\bc(\frak{t},-\alpha).$$
with respect to $\bc Y_1$. We have such a restricted root space decomposition only when $(\frak{s}_\bc,\frak{s}_\bc\cap\frak{h}'_\bc)\cong (\frak{so}(n+2,\bc), \frak{so}(n+1,\bc))$ where $n=\dim_\bc \frak{m}_\bc(\frak{t},\alpha)$.

By the classification in \cite{Be}, we have
$$(\frak{s},\frak{s}^\tau)\cong (\frak{so}(p,q),\frak{so}(p-1,q))\qquad(q=n+2-p)$$
with some $p=1,\ldots,n+1$ since $\frak{s}$ is noncompact. (This can be also deduced from \cite{Mc} by considering two commuting involutions $\tau$ and $\theta$ on the compact real form $(\frak{s}\cap\frak{k})\oplus i(\frak{s}\cap\frak{m})\cong \frak{so}(n+2)$.) If $p=1$, then we have $\frak{s}\cap\frak{k}= \frak{s}\cap\frak{h}'$, a contradiction to $\br iY_1=\frak{s}\cap\frak{t}\subset \frak{k}\cap\frak{q}'$. Hence $p\ge 2$. Taking $\frak{s}\cap\frak{t}=\br(E_{12}-E_{21})\subset \frak{s}\cap\frak{q}'$, we can compute
$$\dim_\bc (\frak{s}_\bc (\frak{t},\alpha)\cap \frak{k}_\bc)=p-2.$$
Since $\frak{s}_\bc (\frak{t},\alpha)\subset \frak{m}_\bc$, we have $p=2$.
\end{proof}

For $Z\in \frak{m}$ we define a parabolic subgroup $P_Z$ of $G$ by
$$P_Z=Z_G(Z)\exp\frak{n}_Z$$
where
$\frak{n}_Z=\bigoplus_{\lambda>0} \{X\in\frak{g}\mid [Z,X]=\lambda X\}$.
For $\alpha\in\Sigma$ let $Y_\alpha$ denote the unique element in $\frak{t}$ such that $\alpha(Y_\alpha)=\pi i/2$ and that $B(Y_\alpha,\frak{t}_\alpha)=\{0\}$ where $\frak{t}_\alpha=\{Y\in\frak{t}\mid \alpha(Y)=0\}$. Put $a_\alpha=\exp Y_\alpha$.

\begin{theorem} \ Suppose that $\frak{m}_\bc(\frak{t},\alpha)\ne 0$ for all the longest roots $\alpha$ in $\Sigma$ and that all the longest roots in $\Sigma$ are mutually conjugate under $W_{K\cap H}(\frak{t})$. Take a longest root $\alpha$ in $\Sigma$ and a nonzero element $Z$ of $(\frak{m}_\bc(\frak{t},\alpha)\oplus \frak{m}_\bc(\frak{t},-\alpha))\cap\frak{h}'$. Then we have$:$

{\rm (i)} \ $H'P_{\pm Z}$ and $Ha_\alpha^{-1}P_{\pm Z}$ are closed in $G$.

{\rm (ii)} \ For any element $x$ in the boundary of $D$, we have
$$xHa_\alpha^{-1}P_Z\cap H'P_Z\ne\phi
\qquad\mbox{or}\qquad
xHa_\alpha^{-1}P_{-Z}\cap H'P_{-Z}\ne\phi.$$

{\rm (iii)} \ Assume moreover that
$\dim\frak{m}_\bc(\frak{t},\alpha)\ge 2$
or that there exists a $t\in T\cap H$ such that $\Ad(t)Z=-Z$. 
Then for any element $x$ in the boundary of $D$, we have
$$xHa_\alpha^{-1}P_Z\cap H'P_Z\ne\phi.$$
\end{theorem}

\begin{proof} (i) \ Since $Z\in\frak{m}\cap\frak{h}'$, it follows that $H'\cap P_Z$ is a parabolic subgroup of $H'$. Hence $H'P_Z/P_Z\cong H'/H'\cap P_Z$ is compact. This implies that $H'P_Z$ is closed in $G$. Replacing $Z$ by $-Z$, we see that $H'P_{-Z}$ is also closed in $G$.

Identify $\frak{s}$ with $\frak{so}(2,n)$ as in Proposition 3.1. Then we have
$$\frak{s}\cap\frak{t}=\br(E_{21}-E_{12})$$
where $E_{jk}$ denote the matrix units. We can assume that $\alpha\in\Sigma$ satisfies
$\alpha: y(E_{21}-E_{12})\mapsto iy.$
Then we have
$Y_\alpha={\pi\over 2}(E_{21}-E_{12})$
and therefore the adjoint action of $a_\alpha=\exp Y_\alpha$ on $\frak{s}$ is equal to that of
$$\bp 0 & -1 & 0 \\ 1 & 0 & 0 \\ 0 & 0 & I_n \ep$$
on $\frak{so}(2,n)$ by the identification $\frak{s}\cong \frak{so}(2,n)$.

On the other hand, $Z$ is of the form
$$\bp 0 & 0 & 0 & \cdots & 0 \\
0 & 0 & x_1 & \cdots & x_n \\
0 & x_1 & 0 & \cdots & 0 \\
\vdots & \vdots & \vdots && \vdots \\
0 & x_n & 0 & \cdots & 0 \ep$$
with some $(x_1,\ldots,x_n)\in\br^n-\{0\}$. Put $Z_0=E_{23}+E_{32}$. Then we see that $Z$ is $Z_{K\cap H}(\frak{t})$-conjugate to $rZ_0$ with some $r\in\br^\times$. (We can take $r>0$ if $n\ge 2$.) Since $P_{cZ}=P_Z$ for $c>0$ and since we consider $\pm Z$, we may assume that
$$Z=Z_0.$$

We see that
$$a_\alpha^{-1}P_{\pm Z}a_\alpha=P_{\pm Z'}$$
where
$$Z'=\Ad(a_\alpha)^{-1}Z=E_{13}+E_{31} \in\frak{s}\cap\frak{m}\cap\frak{q}'=\frak{s}\cap \frak{m}\cap \frak{h}.$$
Hence $Ha_\alpha^{-1}P_{\pm Z}=HP_{\pm Z'}a_\alpha^{-1}$ are also closed in $G$ by the same reason as for $H'P_{\pm Z}$.

(ii) \ If $h\in H$ and $h'\in H'$, then we have
$$(h'xh)Ha_\alpha^{-1}P_{\pm Z}\cap H'P_{\pm Z} =h'(xHa_\alpha^{-1}P_{\pm Z}\cap H'P_{\pm Z}).$$
So we may replace $x$ by any element in the double coset $H'xH$.

By Remark 2.9 (iii) we have $\partial D=(\bigcup_{\Delta\in\mathcal{B}_\alpha} M(\Delta,1))^{cl}$. First assume that $x\in M(\Delta,1)$ for some $\Delta\in\mathcal{B}_\alpha$. Then we have $x\in H'(\exp X_n)aH$ for some $a=\exp Y\in\exp\Delta$ and a nilpotent element $X_n$ in $\frak{q}'\cap\Ad(a)\frak{q}$. By the above remark we may assume that
$$x=(\exp X_n)a.$$
By Lemma 2.7 we have
$$\frak{g}_\bc^{f_a} =\frak{z}_{\frak{k}_\bc}(\frak{t})\oplus \frak{m}_\bc(\frak{t},\alpha) \oplus \frak{m}_\bc(\frak{t},-\alpha).$$

Since $\frak{s}_\bc$ is generated by $\frak{m}_\bc(\frak{t},\alpha) \oplus \frak{m}_\bc(\frak{t},-\alpha)$, it is contained in $\frak{g}_\bc^{f_a}$. Since
$$[\frak{z}_{\frak{k}_\bc}(\frak{t}), \frak{s}_\bc]\subset \frak{s}_\bc,$$
$\frak{s}_\bc$ is an ideal of $\frak{g}_\bc^{f_a}$. Hence $\frak{s}=\frak{s}_\bc\cap\frak{g}$ is an ideal of $\frak{g}^{f_a}$. Let $\frak{s}^\perp$ denote the orthogonal complement of $\frak{s}$ in $\frak{g}^{f_a}$ with respect to the Killing form on $\frak{g}$. Then $\frak{s}^\perp$ is an ideal of $\frak{g}^{f_a}$ contained in $\frak{z}_\frak{k}(\frak{t})$. (Note that $\frak{s}$ may be smaller than $\frak{s}_\Delta$ in Section 2.)

Since $\frak{g}^{f_a}=\frak{s}\oplus \frak{s}^\perp$ and $\frak{s}^\perp\subset \frak{k}$, the nilpotent element $X_n$ is contained in $\frak{s}\cap\frak{q}'$. By the identification in Proposition 3.1, $X_n$ is of the form
$$\bp 0 & -x_1 & x_2 & \cdots & x_{n+1} \\
x_1 &&&& \\
x_2 &&&& \\
\vdots &&& 0 & \\
x_{n+1} &&&& \ep$$
with $x_j\in\br$ such that $x_1^2=x_2^2+\cdots+x_{n+1}^2$. Since we can consider any element in the double coset $H'xH$, we can replace $X_n$ by an $H'\cap aHa^{-1}$-conjugate. Since $\frak{h}'\cap\Ad(a)\frak{h} =\frak{g}^{f_a}\cap\frak{h}'\supset \frak{s}\cap\frak{h}'\cong \frak{so}(1,n)$, we may replace $X_n$ by some $SO(1,n)_0$-conjugate by the identification in Proposition 3.1. So we may assume
$$X_n=\pm(E_{21}-E_{12}+\varepsilon E_{31}+\varepsilon E_{13})$$
where $\varepsilon=\pm 1$. (We may put $\varepsilon=1$ if $n\ge 2$.)
By computation we see that
$$[Z,X_n]=\varepsilon X_n.$$

Since $\alpha(Y)=\alpha(Y_\alpha)=\pi i/2$, we can write $Y=Y_\alpha+Y'$ with $Y'\in\frak{t}\cap\frak{s}^\perp$. Hence we can write
$$a=a_\alpha a'$$
with $a'=\exp Y'$. 

Since $X_n\in\frak{n}_{\varepsilon Z}$ and since $Y'\in\frak{s}^\perp \subset \frak{z}_\frak{g}(Z)$, it follows that
$$xHa_\alpha^{-1}P_{\varepsilon Z}\supset xa_\alpha^{-1}P_{\varepsilon Z}=(\exp X_n)a'P_{\varepsilon Z}=P_{\varepsilon Z}.$$
Hence we have
\begin{equation}
xHa_\alpha^{-1}P_{\varepsilon Z}\cap H'P_{\varepsilon Z}\ne \phi \tag{3.4}
\end{equation}
for $\varepsilon=1$ or $-1$.

Finally let $y$ be an arbitrary element in the boundary of $D$. If $yHa_\alpha^{-1}P_{\pm Z}\cap H'P_{\pm Z}=\phi$, then $xHa_\alpha^{-1}P_{\pm Z}\cap H'P_{\pm Z}=\phi$ for all the elements $x$ in a neighborhood $U_y$ of $y$ because $yHa_\alpha^{-1}P_Z/P_Z$ and $H'P_Z/P_Z$ are compact. But this contradicts (3.4) because $\partial D=(\bigcup_{\Delta\in\mathcal{B}_\alpha} M(\Delta,1))^{cl}$.

The assertion (iii) is also proved in (ii).
\end{proof}

\section{Proof of Theorem 1.3}

In this section, we will prove Theorem 1.3. So we assume that $G$ is a complex Lie group and that $H'\ (=G_\br)$ is a real form of $G$.

Since $K$ is a compact real form of $G$, we have
$$\dim_\bc \frak{m}_\bc(\frak{t},\alpha)
= \dim_\bc \frak{k}_\bc(\frak{t},\alpha)
= \dim_\bc \frak{g}(\frak{t},\alpha)$$
for all $\alpha\in\Sigma(\frak{t})$. Hence we can identify $\Sigma(\frak{m}_\bc,\frak{t})$ with the usual restricted root system $\Sigma=\Sigma(\frak{t})$ of $H'$. It is also known that all the longest roots in $\Sigma$ are mutually conjugate under $W_{K\cap H}(\frak{t})$ for simple $H'$. So the conditions in Theorem 3.2 are satisfied.

Here we give the well-known list of simple real Lie algebras $\frak{h}'$ and the multiplicities $n$ of the longest restricted roots. (Assume $p\le q$.)

\bigskip
\centerline{
\vbox{\offinterlineskip
\hrule
\halign{&\vrule#&\strut\ $\hfil#\hfil$\ \cr
height2pt&\omit&&\omit&&\omit&&\omit&&\omit&\cr
& {\rm type} && \frak{h}' && \Sigma && n && &\cr
height2pt&\omit&&\omit&&\omit&&\omit&&\omit&\cr
\noalign{\hrule}
height4pt&\omit&&\omit&&\omit&&\omit&&\omit&\cr
& {\rm AI} && \frak{sl}(\ell,\br) && {\rm A}_{\ell-1} && 1 && &\cr
height4pt&\omit&&\omit&&\omit&&\omit&&\omit&\cr
& {\rm AII} && \frak{gl}(\ell,\bh)/\br && {\rm A}_{\ell-1} && 4 && &\cr
height4pt&\omit&&\omit&&\omit&&\omit&&\omit&\cr
& {\rm AIII} && \frak{su}(p,q) && {\rm BC}_p\mbox{ or }{\rm C}_p && 1 && {\rm Hermitian} &\cr
height4pt&\omit&&\omit&&\omit&&\omit&&\omit&\cr
& {\rm BDI} && \frak{so}(p,q)\ (p\ge 2) && {\rm B}_p\mbox{ or }{\rm D}_p && 1 && \mbox{Hermitian if $p=2$} &\cr
height4pt&\omit&&\omit&&\omit&&\omit&&\omit&\cr
& {\rm BDI} && \frak{so}(1,q) && {\rm B}_1 && q-1 && &\cr
height4pt&\omit&&\omit&&\omit&&\omit&&\omit&\cr
& {\rm CI} && \frak{sp}(\ell,\br) && {\rm C}_\ell && 1 && {\rm Hermitian} &\cr
height4pt&\omit&&\omit&&\omit&&\omit&&\omit&\cr
& {\rm CII} && \frak{sp}(p,q) && {\rm BC}_p\mbox{ or }{\rm C}_p && 3 && &\cr
height4pt&\omit&&\omit&&\omit&&\omit&&\omit&\cr
& {\rm DIII} && \frak{so}^*(4\ell) && {\rm C}_\ell && 1 && {\rm Hermitian} &\cr
height4pt&\omit&&\omit&&\omit&&\omit&&\omit&\cr
& {\rm DIII} && \frak{so}^*(4\ell+2) && {\rm BC}_\ell && 1 && {\rm Hermitian} &\cr
height4pt&\omit&&\omit&&\omit&&\omit&&\omit&\cr
& {\rm EI} && && {\rm E}_6 && 1 && &\cr
height4pt&\omit&&\omit&&\omit&&\omit&&\omit&\cr
& {\rm EII} && && {\rm F}_4 && 1 && &\cr
height4pt&\omit&&\omit&&\omit&&\omit&&\omit&\cr
& {\rm EIII} && && {\rm BC}_2 && 1 && {\rm Hermitian} &\cr
height4pt&\omit&&\omit&&\omit&&\omit&&\omit&\cr
& {\rm EIV} && && {\rm A}_2 && 8 && &\cr
height4pt&\omit&&\omit&&\omit&&\omit&&\omit&\cr
& {\rm EV} && && {\rm E}_7 && 1 && &\cr
height4pt&\omit&&\omit&&\omit&&\omit&&\omit&\cr
& {\rm EVI} && && {\rm F}_4 && 1 && &\cr
height4pt&\omit&&\omit&&\omit&&\omit&&\omit&\cr
& {\rm EVII} && && {\rm C}_3 && 1 && {\rm Hermitian} &\cr
height4pt&\omit&&\omit&&\omit&&\omit&&\omit&\cr
& {\rm EVIII} && && {\rm E}_8 && 1 && &\cr
height4pt&\omit&&\omit&&\omit&&\omit&&\omit&\cr
& {\rm EIX} && && {\rm F}_4 && 1 && &\cr
height4pt&\omit&&\omit&&\omit&&\omit&&\omit&\cr
& {\rm FI} && && {\rm F}_4 && 1 && &\cr
height4pt&\omit&&\omit&&\omit&&\omit&&\omit&\cr
& {\rm FII} && && {\rm BC}_1 && 7 && &\cr
height4pt&\omit&&\omit&&\omit&&\omit&&\omit&\cr
& {\rm G} && && {\rm G}_2 && 1 && &\cr
height4pt&\omit&&\omit&&\omit&&\omit&&\omit&\cr
& \mbox{complex cases} && && && 2 && &\cr
height4pt&\omit&&\omit&&\omit&&\omit&&\omit&\cr}
\hrule}
}

\bigskip
In the following arguments we consider the complex structure only inside $\frak{g}$. It means that we do not consider the ``complexification'' of the complex Lie algebra $\frak{g}$ in order to avoid confusion.

We define a maximal abelian subspace $\frak{a}$ of $\frak{m}\cap\frak{q}$ by
$$\frak{a}=\br Z\oplus i\frak{t}_\alpha$$
where $\frak{t}_\alpha=\{Y\in\frak{t}\mid \alpha(Y)=0\}$. Take a positive system $\Sigma(\frak{a})^+$ of the restricted root system $\Sigma(\frak{a})=\Sigma(\frak{g},\frak{a})$ so that
$$\beta(Z)>0\mbox{ for all }\beta\in\Sigma(\frak{a})^+.$$
Then the pair $(\frak{a},\Sigma(\frak{a})^+)$ defines a parabolic subgroup
$$P=P(\frak{a},\Sigma(\frak{a})^+)$$
of $G$ contained in $P_Z$.

Let $\frak{j}$ be a maximal abelian subspace of $\frak{m}$ containing $\frak{a}$. 
Let $\Sigma(\frak{j})^+$ be a positive system of the root system $\Sigma(\frak{j})=\Sigma(\frak{g},\frak{j})$ which is compatible with $\Sigma(\frak{a})^+$. Then the pair $(\frak{j},\Sigma(\frak{j})^+)$ defines a Borel subgroup
$$B=B(\frak{j},\Sigma(\frak{j})^+)$$
of $G$ contained in $P$. Since $\frak{a}=\frak{j}\cap\frak{q}$ is maximal abelian in $\frak{m}\cap\frak{q}$ and since $\Sigma(\frak{j})^+$ is $\tau$-compatible, it follows that $H'B$ is closed in $G$ and that $HB$ is open in $G$ (\cite{M1} Proposition 1 and Proposition 2). Since $H$ is a complex symmetric subgroup of $G$, $HB$ is the unique open $H$-$B$ double coset in $G$ and therefore $H'B$ is the unique closed $H'$-$B$ double coset in $G$.

\begin{remark} \ Since Theorem 1.3 concerns orbits on the flag manifold $X=G/B$, we may replace $B$ with any conjugate $pBp^{-1}$ for $p\in P$.
\end{remark}

Let $x$ be an element in the boundary of $D$. By Theorem 3.2 we have
\begin{equation}
xHa_\alpha^{-1}P_Z\cap H'P_Z\ne\phi \tag{4.1}
\end{equation}
or
\begin{equation}
xHa_\alpha^{-1}P_{-Z}\cap H'P_{-Z}\ne\phi. \tag{4.2}
\end{equation}
If $G_\br$ is of non-Hermitian type, then the condition in Theorem 3.2 (iii) is satisfied by Lemma 7.2 in the appendix. Hence we have (4.1) in this case.

Assume (4.1). Then we will show
\begin{equation}
x(Ha_\alpha^{-1}pB)^{cl}\cap H'B\ne\phi \tag{4.3}
\end{equation}
for any $p\in P$. It follows from Lemma 7.1 (i) that $HP=HB$. So we have $H'P=H'B$ (\cite{M2}). Hence (4.3) is equivalent to
\begin{equation}
x(Ha_\alpha^{-1}P)^{cl}\cap H'P\ne\phi \notag
\end{equation}
because $(Ha_\alpha^{-1}P)^{cl}=(Ha_\alpha^{-1}pB)^{cl}P$.
We have only to show
$$(Ha_\alpha^{-1}P)^{cl}=Ha_\alpha^{-1}P_Z$$
which is equivalent to
$(HP')^{cl}=HP_{Z'}$
where $P'=a_\alpha^{-1}Pa_\alpha$. Put
\begin{equation}
\frak{a}'=\Ad(a_\alpha^{-1})\frak{a}=\br Z'\oplus i\frak{t}_\alpha. \tag{4.4}
\end{equation}
Then the Lie algebra $\frak{p}'$ of $P'$ is defined by the pair $(\frak{a}',\Sigma(\frak{a}')^+)$ where $\Sigma(\frak{a}')^+=\{\alpha\circ\Ad(a_\alpha)\in\Sigma(\frak{a}')\mid \alpha\in\Sigma(\frak{a})^+\}$. Since $H$ and $P'$ are complex subgroups of $G$, we have only to show the equality
\begin{equation}
\frak{h}+\frak{p}'=\frak{h}+\frak{p}_{Z'}. \tag{4.5}
\end{equation}
We can write
$$\frak{p}_{Z'}=\frak{p}'\oplus \bigoplus_{\beta\in\Sigma(\frak{a}')^+, \beta(Z')=0} \frak{g}(\frak{a}',-\beta).$$
If $\beta\in\Sigma(\frak{a}')^+$ satisfies $\beta(Z')=0$, then it follows from (4.4) that $\sigma\beta=-\beta$ and hence
$\sigma\frak{g}(\frak{a}',-\beta)=\frak{g}(\frak{a}',\beta)$.
So we have
$\frak{p}_{Z'}\subset\frak{p}'+\sigma\frak{p}'\subset \frak{h}+\frak{p}'$.
(4.5) is clear by this inclusion. Thus we have proved (4.3).

Next assume (4.2). Since $\Ad(a_\alpha^2)Z=-Z$, we have $P_{-Z}=a_\alpha^2P_Za_\alpha^{-2}$. Hence
$$xHa_\alpha^{-1}a_\alpha^2P_Za_\alpha^{-2}\cap H'a_\alpha^2P_Za_\alpha^{-2}\ne\phi$$
and therefore
$$xHa_\alpha P_Z\cap H'a_\alpha^2P_Z\ne\phi.$$
Since $H'P_Z$ and $H'a_\alpha^2P_Z=H'P_{-Z}a_\alpha^2$ are closed in $G$ and since there is only one closed $H'$-$P_Z$ double coset in $G$, we have $H'P_Z=H'a_\alpha^2P_Z$. Hence
$$xHa_\alpha P_Z\cap H'P_Z\ne\phi.$$
By the same argument as for (4.1) $\Longrightarrow$ (4.3), we get
$$x(Ha_\alpha p'B)^{cl}\cap H'B\ne\phi$$
for an arbitrary $p'\in P$. Thus we have proved (ii) for $\widetilde{S}=Ha_\alpha^{-1}pB\cup Ha_\alpha p'B$ or $\widetilde{S}=Ha_\alpha^{-1}pB$ for arbitrary $p$ and $p'$ in $P$.

\begin{remark} \ 
%In the above argument for (ii), we need no assumption on the choice of $B$ in $P$. 
In the following proof of (iii), we must choose $p\in P$ (resp. $p'\in P$) so that $Ha_\alpha^{-1}pB$ (resp. $Ha_\alpha p'B$) is as small as possible. Of  course we may also replace $B$ with $pBp^{-1}$ for some $p\in P$. The following considerations on orbit structure are based on many examples computed in \cite{MO} Section 4.
\end{remark}

Now we will prove (iii). First consider the case where $n$ is odd. Then in view of Proposition 3.1 and the choice of $Z$ in the proof of Theorem 3.2, we can take a maximal abelian subalgebra $\frak{t}_\frak{s}$ of $\frak{s}\cap\frak{k}\cap\frak{h} \ (\subset\frak{z}_\frak{g}(\frak{t}))$ so that $[Z,\frak{t}_\frak{s}]=\{0\}$. We may choose $\frak{j}$ so that it contains $\frak{a}\oplus i\frak{t}_\frak{s}$. We see that $\Ad(a_\alpha^{-2})$ defines the reflection with respect to $Z$ on $\frak{j}\cap\frak{s}_\bc$. Since $\Ad(a_\alpha^{-2})$ acts trivially on $\frak{s}^\perp_\bc$, it acts on $\frak{j}$ as the reflection $w_Z$ with respect to $Z$.

Let $\Psi$ denote the set of simple roots in $\Sigma(\frak{j})^+$. For each subset $\Theta$ of $\Psi$ there corresponds a parabolic subgroup
$$P_\Theta=BW_\Theta B$$
of $G$ where $W_\Theta$ is the subgroup of $W$ generated by $\{w_\beta\mid \beta\in\Theta\}$. We will show that
\begin{equation}
\Theta\ne\Psi\Longrightarrow H'a_\alpha^{-1}B\cap H'P_\Theta=\phi \notag
\end{equation}
because the right hand side is equivalent to $Ha_\alpha^{-1}B\cap HP_\Theta=\phi$ (\cite{M2}).

We can study $H'$-$B$ double cosets by using the Bruhat decomposition as in \cite{Sp}. By the map
$$y\mapsto \tau(y)^{-1}y,$$
the double coset $H'yB$ is mapped into $\tau(B)\tau(y)^{-1}yB$. Since $\Sigma(\frak{j})^+$ is $\tau$-compatible, we can write
$$\tau(B)=w_c^{-1}Bw_c$$
where $w_c$ is the longest element in the Weyl group $W_c$ generated by the reflections with respect to the compact roots in $\Sigma(\frak{j})$. Hence by the map
$$y\mapsto w_c\tau(y)^{-1}y,$$
$H'yB$ is mapped into
$$w_c\tau(B)\tau(y)^{-1}yB=Bw_c\tau(y)^{-1}yB.$$

By this map $H'a_\alpha^{-1}B$ is mapped into $Bw_ca_\alpha^{-2}B$. On the other hand, $H'P_\Theta$ is mapped into
$$w_c\tau(P_\Theta)w_c^{-1}w_cP_\Theta=P_{\Theta'}w_cP_\Theta.$$
Here we write $w_c\tau(P_\Theta)w_c^{-1}=P_{\Theta'}$ with some $\Theta'\subset\Psi$ since $w_c\tau(P_\Theta)w_c^{-1}$ contains $w_c\tau(B)w_c^{-1}=B$. Since $w_ca_\alpha^{-2}$ normalizes $\frak{j}$, it is considered an element of $W$. So we have only to prove that
\begin{equation}
w_ca_\alpha^{-2}\notin W_{\Theta'}W_cW_\Theta. \tag{4.6}
\end{equation}

Let $Z_\Theta$ (resp. $Z_{\Theta'}$) be a nonzero element in $\frak{j}$ which is dominant with respect to $\Sigma(\frak{j})^+$ and $W_\Theta$-invariant (resp. $W_{\Theta'}$-invariant). Then we have
\begin{equation}
B(Z,Z_\Theta)>0\mbox{ (resp. }B(Z,Z_{\Theta'})>0) \tag{4.7}
\end{equation}
with respect to the Killing form $B(\ ,\ )$ as follows. $Z$ is identified by $B(\ ,\ )$ with a dominant root $\widetilde{\alpha}$ in $\Sigma(\frak{j})^+$ (modulo positive constant). It is known that
$\widetilde{\alpha}=\sum_{\beta\in\Psi} m_\beta \beta$
with $m_\beta>0$ for all $\beta\in\Psi$. Since
$\beta(Z_\Theta)\ge 0\mbox{ for all }\beta\in\Psi$
and since
$\beta(Z_\Theta)>0\mbox{ for some }\beta\in\Psi,$
it follows that
$$\widetilde{\alpha}(Z_\Theta)=\sum_{\beta\in\Psi}m_\beta \beta(Z_\Theta)>0.$$
Since $Z\in\frak{a}$, we have
\begin{equation}
w_c(Z)=Z. \tag{4.8}
\end{equation}
Since $Z_\Theta$ and $Z_{\Theta'}$ are dominant for $\Sigma(\frak{j})^+$, we have
\begin{equation}
B(Z_{\Theta'},\Ad(w_c)Z_\Theta)\le B(Z_{\Theta'},\Ad(w_1)Z_\Theta) \tag{4.9}
\end{equation}
for any $w_1\in W_c$.

It follows from (4,7), (4.8) and (4.9) that
\begin{align*}
B(Z_{\Theta'},\Ad(w_c)\Ad(a_\alpha^{-2})Z_\Theta) 
& =B(Z_{\Theta'},\Ad(w_c)(Z_\Theta-{2B(Z,Z_\Theta)\over B(Z,Z)}Z)) \\
& =B(Z_{\Theta'},\Ad(w_c)Z_\Theta)-{2B(Z,Z_\Theta)B(Z,Z_{\Theta'})\over B(Z,Z)} \\
& <B(Z_{\Theta'},\Ad(w_1)Z_\Theta) \\
& =B(Z_{\Theta'},\Ad(w_{\Theta'}w_1w_\Theta)Z_\Theta)
\end{align*}
for all $w_\Theta\in W_\Theta,\ w_1\in W_c$ and $w_{\Theta'}\in W_{\Theta'}$. Thus we have (4.6).

We can also prove that
\begin{equation}
\Theta\ne\Psi\Longrightarrow H'a_\alpha B\cap H'P_\Theta=\phi \notag
\end{equation}
by the same argument. Thus we have proved (iii) when $n$ is odd.

\section{Proof of ${\rm (iii)}$ for even $n$}

In this section we will prove (iii) of Theorem 1.3 when $n$ is even. In this case  there are no real roots in $\Sigma(\frak{j})$ and hence
\begin{equation}
\mbox{there is only one $K\cap H$-conjugacy class of $\sigma$-stable maximal abelian subspace of $\frak{m}$} \tag{5.1}
\end{equation}
(\cite{M1} Theorem 2, \cite{Su} Theorem 6).

By the classification, there are four cases
$$\mbox{$H'\ (=G_\br)$ is complex,\quad AII,\quad $\frak{so}(1,q)$\ ($q$ is odd) and\quad EIV.}$$

\bigskip
\noindent {\bf 5.1. Complex cases}

\bigskip
Let $G_1$ be a complex simple Lie group. Let $K_1$ be a compact real form of $G_1$ and $\theta_1$ the conjugation of $G_1$ with respect to $K_1$. Then $\frak{g}_1=\frak{k}_1\oplus\frak{m}_1$ is a Cartan decomposition of $\frak{g}_1$ where $\frak{m}_1=i\frak{k_1}$.

Put $G=G_1\times G_1,\ H=\{(g,g)\mid g\in G_1\}$ and $K=K_1\times K_1$. Then
$$H'=\{(g,\theta_1(g))\mid g\in G_1\}.$$
Let $\frak{a}=\br Z\oplus i\frak{t}_\alpha\subset \frak{m}\cap\frak{q}$ be as in Section 4. Since
$\frak{m}\cap\frak{q}=\{(X,-X)\mid X\in\frak{m}_1\},$
we can write
$\frak{a}=\{(X,-X)\mid X\in\frak{j}_1\}$
with some maximal abelian subspace $\frak{j}_1$ of $\frak{m}_1$. 

Write
$Z=(Z_1,-Z_1),\ a_\alpha=(a,a^{-1})=(\exp Y,\exp(-Y)),\ i\frak{t}_\alpha=\{(X,-X)\mid X\in(\frak{j}_1)_\alpha\}$
and
$$Z'=\Ad(a^{-1}_\alpha)Z=(\Ad(a^{-1})Z_1,-\Ad(a)Z_1) =(Z_2,Z_2)\in\frak{m}\cap\frak{h}$$
with $Z_1\in\frak{j}_1,\ Y\in\frak{k}_1$ and $Z_2\in\frak{m}_1$. Then we have
$$\Ad(a^2)Z_1=-\Ad(a)Z_2=-Z_1\mand \Ad(a^2)X=X\mbox{ for }X\in(\frak{j}_1)_\alpha.$$
Hence $\Ad(a^2)|_{\frak{j}_1}$ is the reflection with respect to $Z_1$.

Let $\Sigma_1$ denote the root system of the pair $(\frak{g}_1,\frak{j}_1)$ and $\Sigma_1^+$ a positive system of $\Sigma_1$ such that $Z_1$ is dominant for $\Sigma_1^+$. Let $B_1$ be the Borel subgroup of $G_1$ for the pair $(\frak{j}_1, \Sigma_1^+)$ and put $B=B_1\times\theta_1(B_1)$. By the map
\begin{equation}
H'(y,z)B\mapsto B_1y^{-1}\theta_1(z)B_1 \tag{5.2}
\end{equation}
the decomposition $H'\backslash G/B$ is identified with the Bruhat decomposition $B_1\backslash G_1/B_1$ of $G_1$. Let $\Psi_1$ denote the set of simple roots in $\Sigma_1^+$. Then every parabolic subgroup $P_\Theta$ of $G$ containing $B$ is written as
$$P_\Theta=P_{\Theta_1}\times\theta_1(P_{\Theta_2})$$
with some subsets $\Theta_1$ and $\Theta_2$ of $\Psi_1$ where $P_{\Theta_i}=B_1W_{\Theta_i}B_1$ is the parabolic subgroup of $G_1$ corresponding to $\Theta_i$.

By the identification (5.2) $H'a_\alpha^{-1}B$ and $H'P_\Theta$ are identified with $B_1a\theta_1(a)B_1=B_1a^2B_1=B_1w_{Z_1}B_1$ and $P_{\Theta_1}P_{\Theta_2}=B_1W_{\Theta_1}W_{\Theta_2}B_1$, respectively. Since we assume that $H'\backslash G/P_\Theta$ is nontrivial, $\Theta_1$ and $\Theta_2$ are not equal to $\Psi_1$. So we have only to show
$$w_{Z_1}\notin W_{\Theta_1}W_{\Theta_2}.$$
But we can prove this by the same argument as in Section 4 since $Z_1$ corresponds to the maximal root in $\Sigma_1^+$. Thus we have proved (iii) when $H'$ is a complex Lie group.

\bigskip
\noindent {\bf 5.2. AII-case}

\bigskip
Since $n=4$, the pair $(\frak{s},\frak{s}^\tau)$ is isomorphic to $(\frak{so}(2,4),\frak{so}(1,4))$ by Proposition 3.1. By this identification, we can take a three-dimensional maximal abelian subspace 
$$\frak{j}_\frak{s}=\br Z\oplus \br Z_1\oplus \br Z_2$$
in $\frak{m}\cap\frak{s}_\bc$ where $Z=E_{23}+E_{32}$ (as in Section 3), $Z_1=E_{14}+E_{41}$ and $Z_2=i(E_{56}-E_{65})$. Take a maximal abelian subspace $\frak{j}_{\frak{s}^\perp}$ of $\frak{m}\cap(\frak{s}^\perp)_\bc=i\frak{s}^\perp$ containing $i\frak{t}_\alpha$. Then $\frak{j}=\frak{j}_\frak{s}\oplus \frak{j}_{\frak{s}^\perp}$ is a maximal abelian subspace of $\frak{m}$ containing $\frak{a}=\br Z\oplus i\frak{t}_\alpha$. By computation we have
$$\Ad(a_\alpha^{-2})Z=-Z,\quad \Ad(a_\alpha^{-2})Z_1=-Z_1,\quad \Ad(a_\alpha^{-2})Z_2=Z_2$$
and $\Ad(a_\alpha^{-2})X=X$ for all $X\in\frak{j}_{\frak{s}^\perp}$.

Let $\Sigma(\frak{a})^+$ be a positive system of $\Sigma(\frak{a})$ such that $Z$ is dominant for $\Sigma(\frak{a})^+$. Let $\Sigma(\frak{j})^+$ be a positive system of $\Sigma(\frak{j})$ which is compatible with $\Sigma(\frak{a})^+$. Let $\Psi=\{\alpha_1,\ldots, \alpha_m\}$ ($m$ is odd) denote the set of simple roots in $\Sigma(\frak{j})^+$ where
$$(\alpha_j,\alpha_k)=\begin{cases} 1 & \text{if $k=j$}, \\
-1/2 & \text{if $|k-j|=1$}, \\
0 & \text{if $|k-j|\ge 2$} \end{cases}$$
as usual. Then $\{\alpha_k\mid k\mbox{ is odd}\}$ is the set of the compact simple roots in $\Psi$. Let $\alpha_{\rm max}=\alpha_1+\cdots +\alpha_m$ denote the maximal root in $\Sigma(\frak{j})^+$. Then the four roots
$$\alpha_{\rm max},\quad \alpha_{\rm max}-\alpha_1,\quad \alpha_{\rm max}-\alpha_m\mand \alpha_{\rm max}-\alpha_1-\alpha_m$$
are mapped onto the maximal root $\widetilde{\alpha}$ in $\Sigma(\frak{a})^+$ by the restriction of roots in $\Sigma(\frak{j})$ to $\frak{a}$. By the Killing form we can identify $Z$ and $Z_1$ with $\widetilde{\alpha}$ and a vector in $\br\alpha_1\oplus \br\alpha_m$, respectively.

We see that
$$P_Z=P_{\Psi-\{\alpha_2,\alpha_{m-1}\}}=P_{\{\alpha_1\}} P_{\{\alpha_m\}}P_{\Theta_0}$$
where $\Theta_0=\{\alpha_3,\alpha_4, \ldots,\alpha_{m-2}\}$. Since $Ha_\alpha^{-1}P_Z$ is closed in $G$ by Theorem 3.2 (i), it follows that
$$Ha_\alpha^{-1}pP_{\Theta_0}$$
is closed in $G$ for some $p\in P_{\{\alpha_1\}} P_{\{\alpha_m\}}$. It follows from (5.1) and Lemma 7.1 that
$Ha_\alpha^{-1}pB=Ha_\alpha^{-1}wB$
for some $w\in W_{\{\alpha_1\}} W_{\{\alpha_m\}}=\{e, w_{\alpha_1}, w_{\alpha_m}, w_{\alpha_1}w_{\alpha_m}\}$. (More precisely, we can see that $Ha_\alpha^{-1}B=Ha_\alpha^{-1} w_{\alpha_1}w_{\alpha_m}B$ and $Ha_\alpha^{-1}w_{\alpha_1}B=Ha_\alpha^{-1}w_{\alpha_m}B$. But we don't need these equalities here.) Replace $B$ by $wBw^{-1}$. Then
$$Ha_\alpha^{-1}P_{\Theta_0}$$
is closed in $G$.

First suppose that $\Theta$ contains $\Theta_0$. If $Ha_\alpha^{-1}B\subset HP_\Theta$, then $HP_\Theta$ contains a closed set $Ha_\alpha^{-1}P_{\Theta_0}$. Hence $HP_\Theta$ is closed in $G$. But $HP_\Theta$ is also open in $G$. So we have
$HP_\Theta=G$
contradicting the assumption $HP_\Theta\ne G$. Thus we have
$$Ha_\alpha^{-1}B\cap HP_\Theta=\phi.$$

So we have only to consider $\Theta=\Psi-\{\alpha_k\}$ with $k=3,4,\ldots,m-2$. Take a $Z_\Theta$ as in Section 4. Since $w_c\tau(P_\Theta)w_c^{-1}=P_\Theta$, we have only to prove (4.6) for $\Theta'=\Theta$. Since
$\alpha_1(Z_\Theta)=\alpha_m(Z_\Theta)=0$, we have
\begin{equation}
w_{Z_1}Z_\Theta=Z_\Theta. \tag{5.3}
\end{equation}
Since $Z$ corresponds to a positive constant multiple of $2\alpha_{\rm max}-\alpha_1-\alpha_m$ and since
$$(2\alpha_{\rm max}-\alpha_1-\alpha_m)(Z_\Theta)= 2\alpha_{\rm max}(Z_\Theta)>0$$
as in Section 4, we have
\begin{equation}
B(Z,Z_\Theta)>0. \tag{5.4}
\end{equation}
By (4.8), (4.9), (5.3) and (5.4) we have
\begin{align*}
B(Z_\Theta,\Ad(w_ca_\alpha^{-2})Z_\Theta) 
& =B(Z_\Theta,\Ad(w_cw_Zw_{Z_1})Z_\Theta) \\
& =B(Z_\Theta,\Ad(w_cw_Z)Z_\Theta) \\
& =B(Z_\Theta,\Ad(w_c)(Z_\Theta-{2B(Z,Z_\Theta)\over B(Z,Z)}Z)) \\
& =B(Z_\Theta,\Ad(w_c)Z_\Theta)-{2B(Z,Z_\Theta)^2\over B(Z,Z)} \\
& <B(Z_\Theta,\Ad(w_1)Z_\Theta) \\
& =B(Z_\Theta,\Ad(w_\Theta w_1w'_\Theta)Z_\Theta)
\end{align*}
for all $w_\Theta,w'_\Theta\in W_\Theta$ and $w_1\in W_c$. Thus we have proved
$$w_ca_\alpha^{-2}\notin W_\Theta W_c W_\Theta.$$

\bigskip
\noindent {\bf 5.3. Cases of $\frak{so}(1,q)$\ ($q$ is odd) and EIV}

\bigskip
In these cases, we have $P_Z=P$. (Remark: We only use this condition. So the following argument is also valid when $H'\ (=G_\br)$ is real rank one.) Since $Ha_\alpha^{-1}P=Ha_\alpha^{-1}P_Z$ is closed in $G$ by Theorem 3.2 (i), there exists a $p\in P$ such that
$Ha_\alpha^{-1}pB$
is closed in $G$. If $Ha_\alpha^{-1}pB\subset HP_\Theta$, then $HP_\Theta$ is closed in $G$. But $HP_\Theta$ is also open in $G$. Hence we have
$HP_\Theta=G$
a contradiction to the assumption $HP_\Theta\ne G$. Thus we have proved
$$Ha_\alpha^{-1}pB\cap HP_\Theta=\phi$$
and we have completed the proof of Theorem 1.3.

\section{Proof of Theorem 1.7}

Applying Lemma 2 in \cite{GM2} recursively, we can find a sequence of simple roots $\alpha_1,\ldots,\alpha_\ell\ (\ell=\dim_\bc G -\dim_\bc \widetilde{S})$ such that
$$\dim_\bc \widetilde{S}P_{\alpha_1}\cdots P_{\alpha_k}= \dim_\bc \widetilde{S}P_{\alpha_1}\cdots P_{\alpha_{k-1}} +1$$
for $k=1,\ldots,\ell$.

\begin{lemma} \ Suppose that $H'$ is of non-Hermitian type.

{\rm (i)} \ If $H'$ is not of {\rm CII}-type or {\rm FII}-type, then
$$(\widetilde{S}')^{cl}=S'_{\rm op}P_{\alpha_\ell}\cdots P_{\alpha_1}.$$

{\rm (ii)} \ If $H'$ is of {\rm CII}-type, then we can take $\alpha_1$ as a compact root and we have
$${S'_1}^{cl}=S'_{\rm op}P_{\alpha_\ell}\cdots P_{\alpha_2}$$
for the dense $H$-$B$ double coset $S_1$ in $\widetilde{S}P_{\alpha_1}$.
\end{lemma}

\begin{proof} 
Let $S_k$ denote the dense $H$-$B$ double coset in $\widetilde{S}P_{\alpha_1}\cdots P_{\alpha_k}$ for $k=0,\ldots,\ell$. Write $S_k=Hx_kB$ with a representative $x_k$ such that $\frak{j}_k=\Ad(x_k)\frak{j}$ is a $\sigma$-stable maximal abelian subspace in $\frak{m}$. Let $\beta_k$ denote the simple root in
$$\Sigma(\frak{j}_k)^+=\{\gamma\circ\Ad(x_k)^{-1}\mid \gamma\in\Sigma(\frak{j})^+\}$$
defined by $\beta_k=\alpha_k\circ\Ad(x_k)^{-1}$. Since $\dim_\bc S_k=\dim_\bc S_{k-1}+1$, it follows from Lemma 7.1 that
$$\dim(\frak{j}_k\cap\frak{q})=\dim(\frak{j}_{k-1}\cap\frak{q})\mand \beta_k\mbox{ is complex}$$
or that
$$\dim(\frak{j}_k\cap\frak{q})=\dim(\frak{j}_{k-1}\cap\frak{q})+1\mand \beta_k\mbox{ is real}.$$
Especially we have
\begin{equation}
\dim(\frak{j}_0\cap\frak{q})\le \dim(\frak{j}_1\cap\frak{q})\le \cdots \le \dim(\frak{j}_\ell\cap \frak{q}). \tag{6.1}
\end{equation}

Suppose that $\beta_k$ is complex. Then it follows from Lemma 7.1 (iv) that
$$S_kP_{\alpha_k}=S_{k-1}\sqcup S_k.$$
Moreover we can take $x_k=x_{k-1}w_{\alpha_k}$ and hence $\frak{j}_k=\frak{j}_{k-1}$. By the duality (\cite{M2}) we have
$S'_kP_{\alpha_k}=S'_{k-1}\sqcup S'_k$
and $S'_k\subset {S'_{k-1}}^{cl}$. Hence
\begin{equation}
{S'_{k-1}}^{cl}=(S'_kP_{\alpha_k})^{cl}={S'_k}^{cl}P_{\alpha_k}. \tag{6.2}
\end{equation}

(i) \ First suppose that $n$ is even. By (5.1) we have
$$\dim(\frak{j}_0\cap\frak{q})=\dim(\frak{j}_1\cap\frak{q})= \cdots =\dim(\frak{j}_\ell\cap \frak{q}).$$
Since $\beta_k$ is complex for all $k=1,\ldots,\ell$, it follows from (6.2) that
$$(\widetilde{S}')^{cl}={S'_0}^{cl}={S'_1}^{cl}P_{\alpha_1}= \cdots ={S'_\ell}^{cl}P_{\alpha_\ell}\cdots P_{\alpha_1} =S'_{\rm op} P_{\alpha_\ell}\cdots P_{\alpha_1}.$$

Next consider the case where $n$ is odd. Since we assume $H'$ is not of {\rm CII}-type or {\rm FII}-type, it follows that $n=1$. Take a $\sigma$-stable maximal abelian subspace $\frak{j}$ of $\frak{m}$ as in Section 4 and Section 7.2. Then $\frak{j}'=\Ad(a_\alpha^{-1})\frak{j}$ is $\sigma$-stable. So we may assume that $x_0=a_\alpha^{-1}$. Since $\frak{j}'\cap \frak{q}=i\frak{t}_\alpha$, it follows from (6.1) that
\begin{equation}
\dim(\frak{j}_0\cap\frak{q})= \cdots =\dim(\frak{j}_{m-1}\cap\frak{q}) =\dim(\frak{j}_m\cap \frak{q})-1= \cdots =\dim(\frak{j}_\ell\cap \frak{q})-1 \notag
\end{equation}
for some $m$. So we may assume that $\frak{j}_0=\cdots =\frak{j}_{m-1}=\frak{j}'$ and that $\frak{j}_m=\cdots =\frak{j}_\ell=\frak{j}$. We see that the root $\widetilde{\alpha}\in\Sigma(\frak{j})$ corresponding to $Z$ is the real root defining $\frak{j}'=\Ad(a_\alpha^{-1})\frak{j}$ from $\frak{j}$ (\cite{Su} Theorem 6). On the other hand $\beta_m$ is also a real root defining the $K\cap H$-conjugacy class of $\sigma$-stable maximal abelian subspace containing $\frak{j}'$. So $\widetilde{\alpha}$ and $\beta_m$ are conjugate under some  $w\in W_{K\cap H}(\frak{j})$. Replacing $x_m$ by $wx_m$, we may assume that $\beta_m=\widetilde{\alpha}$. Consider the root $\gamma=\alpha_m\circ \Ad(x_{m-1})^{-1}=\beta_m\circ \Ad(x_mx_{m-1}^{-1})$ in $\Sigma(\frak{j}_{m-1})$. This root $\gamma$ is the noncompact simple root in $\Sigma(\frak{j}_{m-1})^+$ corresponding to $Z'$.

Since $H'$ is of non-Hermitian type, it follows from Lemma 7.2 that $Hw_\gamma x_{m-1}B=Hx_{m-1}B$. Hence we have
$$S_{m-1}P_{\alpha_m}=Hx_{m-1}B\sqcup Hc_{\gamma}x_{m-1}B=S_{m-1}\sqcup S_m$$
by Lemma 7.1 (ii). So we have
$S'_mP_{\alpha_m}=S'_m\sqcup S'_{m-1}$
and $S'_m\subset {S'_{m-1}}^{cl}$. Hence
\begin{equation}
{S'_{m-1}}^{cl}=(S'_mP_{\alpha_m})^{cl}={S'_m}^{cl}P_{\alpha_m}. \tag{6.3}
\end{equation}
By (6.2) and (6.3), we get
$$(\widetilde{S}')^{cl}={S'_0}^{cl}={S'_1}^{cl}P_{\alpha_1}= \cdots ={S'_\ell}^{cl}P_{\alpha_\ell}\cdots P_{\alpha_1} =S'_{\rm op} P_{\alpha_\ell}\cdots P_{\alpha_1}.$$

(ii) \ Take $\frak{j}, \ \Sigma(\frak{j})^+$ and $\Psi$ (the set of the simple roots in $\Sigma(\frak{j})^+$) as in Section 4. Then $\frak{j}_0=\Ad(a_\alpha^{-1})\frak{j}$ is $\sigma$-stable and $\dim(\frak{j}_0\cap\frak{q})=\dim(\frak{j}\cap\frak{q})-1$. Write $\Psi=\{e_1-e_2, e_2-e_3, \ldots, e_{r-1}-e_r, 2e_r\}$ as usual by using an orthonormal basis $\{e_1,\ldots, e_r\}$ of $\frak{j}^*$. Then the compact simple roots are
$$e_1-e_2,\ e_3-e_4,\ \ldots,\ e_{2s-1}-e_{2s},\ e_{2s+1}-e_{2s+2},\ e_{2s+2}-e_{2s+3},\ \ldots,\ ,\ e_{r-1}-e_{r},\ 2e_r$$
where $s\ (\le r/2)$ is the real rank of $H'=Sp(s,r-s)$. The dominant root $\widetilde{\alpha}$ corresponding to $Z$ is
$$\widetilde{\alpha}=e_1+e_2.$$
The root $\beta=\widetilde{\alpha}\circ\Ad(a_\alpha)\in \Sigma(\frak{j}_0)$ is noncompact. Suppose that $(e_1-e_2)\circ\Ad(a_\alpha)$ is compact. Then $2e_1\circ\Ad(a_\alpha)$ and $2e_2\circ\Ad(a_\alpha)$ are noncompact. Since these two roots are strongly orthogonal, we can construct a $\sigma$-stable maximal abelian subspace $\frak{j}'$ of $\frak{m}$ such that
$$\dim(\frak{j}'\cap\frak{q})=\dim(\frak{j}_0\cap \frak{q})+2 =\dim(\frak{j}\cap\frak{q})+1.$$
But this contradicts that $\frak{j}\cap\frak{q}=\frak{a}$ is maximal abelian in $\frak{m}\cap\frak{q}$. Thus we have proved that $(e_1-e_2)\circ\Ad(a_\alpha)$ is a noncompact root of $\Sigma(\frak{j}_0)$.

Put $\alpha_1=e_1-e_2$. Then the dense $H$-$B$ double coset $S_1$ in $\widetilde{S}P_{\alpha_1}$ is written as $S_1=Hx_1B$ where $\frak{j}_1=\Ad(x_1)\frak{j}$ satisfies
$$\dim(\frak{j}_1\cap\frak{q})=\dim(\frak{j}_0\cap\frak{q})+1 =\dim(\frak{j}\cap\frak{q}).$$
Take $\alpha_2,\ldots,\alpha_\ell$ so that
$$\dim_\bc \widetilde{S}P_{\alpha_1}\cdots P_{\alpha_k}= \dim_\bc \widetilde{S}P_{\alpha_1}\cdots P_{\alpha_{k-1}} +1$$
for $k=2,\ldots,\ell$. Then we have
$$\dim(\frak{j}_1\cap\frak{q})=\dim(\frak{j}_2\cap\frak{q})= \cdots =\dim(\frak{j}_\ell\cap \frak{q})$$
by (6.1). Hence $\beta_2,\ldots,\beta_\ell$ are complex roots and so we have
$${S'_1}^{cl}={S'_2}^{cl}P_{\alpha_2}= \cdots ={S'_\ell}^{cl}P_{\alpha_\ell}\cdots P_{\alpha_2} =S'_{\rm op} P_{\alpha_\ell}\cdots P_{\alpha_2}$$
by (6.2).
\end{proof}

\bigskip
\noindent {\it Proof of Theorem 1.7}. \ We will show that
$\widetilde{S}\subset T_j=S_j^{cl}$
for every $j\in J$. By \cite{GM2} Lemma 2, there exists a simple root $\beta$ such that $S_jP_\beta \supset S_{\rm op}$. Hence $S_j\subset S_{\rm op}P_\beta$ and therefore
$S'_j\subset S'_{\rm op}P_\beta.$
Especially $\beta$ is not a compact root. 

If $H'$ is of FII-type, then $\Psi(\frak{j})$ is ${\rm F}_4$-type and it consists of three compact roots and one complex root $\beta$. Hence $J$ consists of only one element $j$ and
$S_{\rm op}P_\beta=S_{\rm op}\sqcup S_j$.
So it is clear that
$$\widetilde{S}\subset T_j$$
because $\widetilde{S}\subset T_j$ for some $j\in J$ by \cite{GM2} Theorem 2.

So we may assume that $H'$ is not of FII-type. It is known that
$$\widetilde{S}\subset S_j^{cl}\Longleftrightarrow (\widetilde{S}')^{cl}\supset S'_j$$
(\cite{M3} Corollary, \cite{MUV} Corollary 1.4). Suppose that $\widetilde{S}\not\subset T_j= S_j^{cl}$. Then $(\widetilde{S}')^{cl}\not\supset S'_j$. Hence
$$\alpha_1,\ldots,\alpha_\ell\ne\beta$$
by Lemma 6.1. (If $H'$ is of CII-type, then $\widetilde{S}\not\subset T_j$ implies $S_1\not\subset T_j$ because $\widetilde{S}\subset S_1^{cl}$. Hence ${S'_1}^{cl}\not\supset S'_j$ and therefore $\alpha_2,\ldots,\alpha_\ell\ne\beta$ by Lemma 6.1 (ii). Since $\alpha_1$ is compact, we have $\alpha_1,\ldots,\alpha_\ell\ne\beta$.) Thus we have
$$\widetilde{S}\subset S_{\rm op}P_\Theta$$
with $\Theta=\Psi-\{\beta\}$. In Remark 1.2 (ii) we see that $HP_\Theta=G$ holds only when $\Theta=\Psi-\{\beta\}$ with some compact root $\beta$. So we have $HP_\Theta\ne G$ and it follows from Theorem 1.3 (iii) that
$$\widetilde{S}\cap S_{\rm op}P_\Theta=\phi$$
a contradiction. Thus we have proved $\widetilde{S}\subset T_j$. \hfill $\square$

\section{Appendix}

In this appendix, we assume $H'$ is a simple Lie group and $G$ is a complexification of $H'$. (So $G,\ H$ and $H'$ are $G_\bc,\ K_\bc$ and $G_\br$, respectively, in Section 1.) In Section 7.2 we will moreover assume that $H'$ is of non-Hermitian type.

\bigskip
\noindent {\bf 7.1 Lemma 5.1 in \cite{V}}

\bigskip
First we will review a lemma due to Vogan \cite{V} which is used frequently in this paper. This lemma is generalized for arbitrary real symmetric pairs $(G,H)$ in \cite{M3} Lemma 3. But we will restrict ourselves to complex symmetric pairs for simplicity.

The full flag manifold of $G$ is identified with the set $\mathcal{F}$ of all the Borel subgroups in $G$. Every $H$-orbit in $\mathcal{F}$ contains a Borel subgroup of the form
\begin{equation}
B=B(\frak{j},\Sigma^+)=Z_G(\frak{j})\exp \frak{n} \tag{7.1}
\end{equation}
with
$\frak{n}=\bigoplus_{\alpha\in\Sigma^+} \frak{g}(\frak{j},\alpha)$
where $\frak{j}$ is a $\sigma$-stable maximal abelian subspace of $\frak{m}$ and $\Sigma^+$ is a positive system of the root system $\Sigma$ of the pair $(\frak{g},\frak{j})$ (\cite{M1} Theorem 1, \cite{R} Theorem 13). ($\frak{j}_\bc$ is a $\sigma$-stable Cartan subalgebra of $\frak{g}$.) By the symmetry of $H$ and $H'$ every $H'$-orbit in $\mathcal{F}$ also contains a Borel subgroup $B$ of the form (7.1). In \cite{M1} (Corollary of Theorem 3) we defined the natural correspondence between $H$-orbits $S$ in $\mathcal{F}$ and $H'$-orbits $S'$ in $\mathcal{F}$ so that $S$ and $S'$ contain the same Borel subgroups $B$ of the form (7.1). (We remark here that the $H'$-orbit structure on the full flag manifold $\mathcal{F}$ was first explicitly studied in \cite{A}.)

Roots in $\Sigma$ are usually classified as follows.

(i) \ If $\sigma(\alpha)=\alpha$ and $\frak{g}(\frak{j},\alpha)\subset \frak{h}$, then $\alpha$ is called a ``compact root''.

(ii) \ If $\sigma(\alpha)=\alpha$ and $\frak{g}(\frak{j},\alpha)\subset \frak{q}$, then $\alpha$ is called a ``noncompact root''.

(iii) \ If $\sigma(\alpha)=-\alpha$, then $\alpha$ is called a ``real root''.

(iv) \ If $\sigma(\alpha)\ne \pm\alpha$, then $\alpha$ is called a ``complex root''.

\noindent For a simple root $\alpha$ of $\Sigma^+$, we can define a parabolic subgroup $P_\alpha$ by
$$P_\alpha=B\sqcup Bw_\alpha B.$$

\begin{lemma} \ $($\cite{V} Lemma 5.1, c.f. also \cite{M3} Lemma 3$)$ \ $HP_\alpha$ is decomposed into $H$-$B$ double cosets as follows.

{\rm (i)} \ If $\alpha$ is compact, then $HP_\alpha=HB$.

{\rm (ii)} \ If $\alpha$ is noncompact, then $HP_\alpha=HB\cup Hw_\alpha B\cup Hc_\alpha B$ and $\dim_\bc HB=\dim_\bc Hw_\alpha B=\dim_\bc Hc_\alpha B -1$. Here $c_\alpha=\exp(X+\theta X)$ with some $X\in\frak{g}(\frak{j},\alpha)$ such that $c_\alpha^2=w_\alpha$. $($Sometimes $HB$ and $Hw_\alpha B$ coincide.$)$

{\rm (iii)} \ If $\alpha$ is real, then $HP_\alpha=HB\cup Hc_\alpha B\cup Hc_\alpha^{-1}B$ and $\dim_\bc Hc_\alpha B=\dim_\bc Hc_\alpha^{-1}B =\dim_\bc HB -1$. Here $c_\alpha=\exp(X+\theta X)$ with some $X\in\frak{g}(\frak{j},\alpha)\cap\frak{q}'$ such that $c_\alpha^2=w_\alpha$. $($Sometimes $Hc_\alpha B$ and $Hc_\alpha^{-1}B$ coincide.$)$

{\rm (iv)} \ If $\alpha$ is complex, then $HP_\alpha=HB\sqcup Hw_\alpha B$. Moreover we have
$$\dim_\bc Hw_\alpha B=\begin{cases} \dim_\bc HB +1 & \text{if $\sigma\alpha\in\Sigma^+$,} \\ \dim_\bc HB -1 & \text{if $\sigma\alpha\notin\Sigma^+$.} \end{cases}$$
\end{lemma}

\bigskip
\noindent {\bf 7.2 A lemma for non-Hermitian cases}

\bigskip
Suppose that $n=\dim_\bc \frak{g}(\frak{t},\alpha)=1$ for a longest root $\alpha$ in $\Sigma(\frak{t})$. Let $\widetilde{\frak{t}}$ be a maximal abelian subalgebra of $\frak{k}$ containing $\frak{t}$. Then the restricted root $\alpha$ is the restriction of the root $\widetilde{\alpha}$ for $\widetilde{\frak{t}}$ such that $\widetilde{\alpha}$ vanishes on $\frak{t}_\frak{h}=\widetilde{\frak{t}}\cap \frak{h}$.

As in Section 4 define a maximal abelian subspace $\frak{a}=\br Z\oplus i\frak{t}_\alpha$ of $\frak{m}\cap\frak{q}$. Then $\frak{j}=\frak{a}\oplus\frak{t}_\frak{h}$ is a maximal abelian subspace of $\frak{m}$. Put $\frak{j}'=\Ad(a_\alpha^{-1})\frak{j}= \br Z'\oplus i\frak{t}_\alpha\oplus\frak{t}_\frak{h}$ where $Z'=\Ad(a_\alpha^{-1})Z$ as in Section 3. Then $\frak{j}'$ is a $\sigma$-stable maximal abelian subspace of $\frak{m}$ such that $\dim (\frak{j}'\cap\frak{q})=\dim (\frak{j}\cap\frak{q})-1$.

\begin{lemma} \ Suppose that $H'$ is of non-Hermitian type. Then there exists a $t\in T\cap H$ such that $e^\alpha(t)=-1$ and that $\Ad(t)|_{\frak{j}'}$ is the reflection with respect to $Z'$. 
\end{lemma}

\begin{proof} Clearly we may assume that $G$ is simply connected. So the compact real form $K$ of $G$ is also simply connected. Then it is known that $K^\sigma$ is connected. So we have
$$K^\sigma=K\cap H.$$
Note that $\Sigma(\frak{t})$ is identified with the restricted root system of the compact symmetric pair $(K,K\cap H)$. Let $Z_\beta$ be the element of $\frak{t}$ defined by
$$Z_\beta={4\pi i\beta\over (\beta,\beta)}.$$
(It means $\gamma(Z_\beta)=4\pi i(\gamma,\beta)/(\beta,\beta)$ for all $\gamma\in\Sigma(\frak{t})$.) Then it is known that the lattice
$$\{Y\in\frak{t}\mid \exp Y=e\}$$
is generated by $\{Z_\beta\mid \beta\in\Sigma(\frak{t})\}$ (c.f. \cite{Mc} Appendix).

Note that $\Sigma(\frak{t})$ is also identified with the restricted root system of $H'$. It is known that $H'$ is of Hermitian type if and only if the following two conditions are satisfied.

(i) \ $\Sigma(\frak{t})$ is C-type or BC-type.

(ii) \ $\dim_\bc \frak{g}(\frak{t},\alpha)=1$ for the longest roots $\alpha$ in $\Sigma(\frak{t})$.

\noindent Since we assume $H'$ is of non-Hermitian type and since we assume (ii), the restricted root system $\Sigma(\frak{t})$ does not satisfy (i). Hence there exists a $\beta\in\Sigma(\frak{t})$ such that \begin{equation}
{(\alpha,\beta)\over (\beta,\beta)}={1\over 2}. \tag{7.2}
\end{equation}

Consider the element
$t=\exp(1/2)Z_\beta$
of $T$. Since
$t\sigma(t)^{-1}=t^2=\exp Z_\beta=e$,
we have $t\in K^\sigma=K\cap H$. By (7.2) we have
$$\alpha\left({1\over 2}Z_\beta\right)={2\pi i(\alpha,\beta)\over (\beta,\beta)}=\pi i$$
and therefore $e^\alpha(t)=e^{\pi i}=-1$.

Since $Z'\in\frak{m}_\bc(\frak{t},\alpha)\oplus \frak{m}_\bc(\frak{t},-\alpha)$, we have $\Ad(t)Z'=-Z'$. (Here we consider the complexification of the complex Lie algebra $\frak{g}$.) On the other hand $\Ad(t)$ acts trivially on $i\frak{t}_\alpha\oplus\frak{t}_\frak{h}$. Hence $\Ad(t)$ acts on $\frak{j}'$ as the reflection with respect to $Z'$.
\end{proof}

\end{document}